\documentclass[a4paper,english,11pt]{article}
\usepackage[utf8]{inputenc}
\usepackage{babel}
\usepackage{longtable}
\usepackage{psboxit}
\usepackage{vmargin,hhline,pifont}
\usepackage{enumerate,fancyhdr}
\usepackage{amsmath}
\usepackage{amsthm}
\usepackage{amsfonts}
\usepackage{amssymb}
\usepackage{chapterbib}
\usepackage[arrow, matrix, curve]{xy}

\newcommand{\ndash}{\nobreakdash-\hspace{0pt}}

\newtheoremstyle{axel}
{}		
{}		
{\itshape}	
{}		
{\bfseries}	
{}		
{0pt\penalty1000}		
{}		
\newtheoremstyle{axelnl}
{}		
{}		
{\itshape}	
{}		
{\bfseries}	
{}		
{\newline}	
{}		
\theoremstyle{axelnl}
\newtheorem{theorem}{Theorem}[section]
\newtheorem{definition}{Definition}[section]
\newtheorem{cor}{Corollary}[section]
\newtheorem{remark}{Remark}[section]
\newtheorem{lemma}{Lemma}[section]

\newtheorem{example}{Example}[section]

\theoremstyle{axel}

\begin{document}
\title{Linear representations of twin cities}
\author{Walter Freyn\thanks{Universit\"at M\"unster, Institut f\"ur Mathematik, Einsteinstra\ss e 62, 48149 M\"unster, walter.freyn@uni-muenster.de}}
\maketitle
\pagestyle{plain}
\pagenumbering{arabic}

\begin{abstract}
For spherical Tits buildings of the classical types there are well-known explicit descriptions as flag complexes. Similarly for affine buildings of the classical types there are explicit constructions in terms of lattices. In this article we generalize the flag complex description to twin cities, a generalization of twin buildings adapted to analytic completions of Kac-Moody groups as they appear for example in Kac-Moody geometry. 
We construct linear representations of twin cities as flag complexes of certain subspaces in Hilbert spaces.

\end{abstract}

\section{Introduction}

A priori Tits buildings are abstract objects associated to groups with a BN\ndash pair, for example Lie groups. They are defined either as simplicial complexes or as metric graphs satisfying some lists of axioms~\cite{AbramenkoBrown08}. From this abstract point of view Lie groups are smooth manifolds with a group structure. Complementary to this abstract approach is the construction of explicit realizations. Explicit realizations of Lie groups are linear representations, that is homomorphisms of a Lie group $G$ into the general linear group $GL(V)$ of a suitable vector space $V$. The corresponding realizations of buildings are flag complexes. A flag complex is a simplicial complexes constructed from a set of flags (chains of subspaces $V_{i_1}\subset\dots\subset  V_{i_k}$ such that the spaces $V_{i}$ are subspaces of $V$) satisfying certain technical requirements. Furthermore we want the representation of the building to be equivariant with respect to the action of the representation of its isometry group. To be more precise let $\mathfrak{B}$ be a Tits building, $G$ its isometry group. Let furthermore $\rho: G\longrightarrow GL(V)$ be a representation of  $G$ and define $\mathcal{SUB}(V)$ to be the space of ``admissible'' subspaces of $V$ and $\mathbb{P}(\mathcal{SUB}(V))$ its power set. The representation $\rho_*:\mathfrak{B}\longrightarrow \mathbb{P}(\mathcal{SUB}(V))$ is a map such that the following diagram commutes.

\begin{displaymath}
\begin{xy}
  \xymatrix{
      \mathfrak{B} \ar[rr]^{g\in G} \ar[d]_{\rho_*}   		&  &   \mathfrak{B} \ar[d]^{\rho_{*}}  \\
      \mathbb{P}(\mathcal{SUB}(V) \ar[rr]^{\rho(g)\in GL(V)})         	&  &   \mathbb{P}(\mathcal{SUB}(V)) 
  }
\end{xy}
\end{displaymath}

In this diagram the action of $\rho(g)$ for $g\in G$ on $\mathbb{P}(\mathcal{SUB}(V))$ is the one induced by the action of $\rho(G)$ on $V$.
Linear representations for the Lie groups of types $A_n$, $B_n$, $C_n$ and $D_n$ and their buildings are well known~\cite{Garrett97}. In the situation of Kac-Moody groups and Kac-Moody algebras, there are linear representations known for affine Kac-Moody algebras and Kac-Moody groups: They can be realized as $2$\ndash dimensional extensions of groups (resp.~algebras) of polynomial maps~\cite{Kac90,Carter05,Tits84}. Similarly there are linear representations of various analytic completions of Kac-Moody algebras that appear in infinite dimensional differential geometry or the theory of integrable systems~\cite{PressleySegal86}. Those representations describe Kac-Moody groups as subgroups of the general linear group of Hilbert\ndash, Banach\ndash\, or Fr\'echet\ndash spaces. 
Representations for affine buildings over local fields are described in~\cite{Garrett97}. A twin building of type $\widetilde{A}_n$ is constructed in~\cite{AbramenkoRonan98}. Lattice chain models for the buildings for the classical groups are a closely related variant~\cite{AbramenkoNebe02}.

Nevertheless twin buildings exist only for the algebraic Kac-Moody groups. Every type of completion destroys the twin BN\ndash pair and hence the action on the associated twin building. Hence there are no twin buildings for most of the Kac-Moody groups and loop groups that appear in infinite dimensional differential geometry. We solved this problem by the introduction of (twin) cities in~\cite{Freyn09,Freyn10a,Freyn10d}. In~\cite{Freyn09,Freyn10a} we called cities ``universal geometric twin buildings'', but changed the denomination, being adverted of the risk of confusion with the universal ``algebraic'' twin building, introduced by Ronan and Tits~\cite{RonanTits94}. In those references our description of the cities are based crucially on the affine Kac-Moody groups acting on them. Affine Kac-Moody groups can be realized as certain torus extensions $\widehat{L}(G)$ of loop groups $L(G)$. Affine Kac-Moody algebras are $2$\ndash dimensional extensions of loop algebras $L(\mathfrak{g})$. Here $G$ denotes a compact simple Lie group, $\mathfrak{g}$ its Lie algebra. Depending on the regularity assumptions on the loops (e.g. of Sobolev class $H^k$) one gets families of completions of the minimal (=algebraic) affine Kac-Moody groups $\widehat{L_{alg}G}$. Following Tits those groups are called analytic Kac-Moody groups. Then our result states that associated to every affine Kac-Moody group $\widehat{L}(G)$ there is a twin city such that the following holds~\cite{Freyn10d}:

\begin{theorem}[Twin cities]
   For each analytic Kac-Moody group $\mathcal{G}$ there exists an associated twin city $\mathfrak{B}= \mathfrak{B}^+ \cup \mathfrak{B^-}$, such that
   \begin{enumerate}[(i)]
     \item Each connected component $\Delta\in \mathfrak{B}$ is an affine building.
     \item Each pair $(\Delta^+, \Delta^-)\in \mathfrak{B}^+\cup
      \mathfrak{B}^-$, consisting of a building $\Delta^+$ in $\mathfrak{B}^+$ (``positive'' building) and a building $\Delta^-$ in $\mathfrak{B}^-$ (``negative'' building),
     is an affine twin building.
     \item $\mathcal{G}$ acts on its twin city $\mathfrak{B}$ by isometries.
     \item ``Small'' twin cities, associated to Kac-Moody groups, defined by stronger regularity conditions, embed into ``big'' twin cities, associated to Kac-Moody groups, defined by weaker regularity conditions.
   \end{enumerate}
\end{theorem}

In this article we construct linear representations for twin cities associated to analytic affine Kac-Moody groups of the types  $\widetilde{A}_n$, $\widetilde{B}_n$, $\widetilde{C}_n$ and $\widetilde{D}_n$. From a structural point of view those affine Kac-Moody groups are closely related to the corresponding finite dimensional Lie algebras of the types $A_n$, $B_n$, $C_n$ and $D_n$: From an algebraic point of view, this can be seen as the Cartan matrix of affine type $\widetilde{X}_n$ for $X\in \{A, B, C, D\}$ can be constructed from the Cartan matrix of finite  type $X_n$ by the addition of an single additional row and column. Thus the root systems and the Weyl groups are closely related: The Weyl group of $\widetilde{X}_n$ is the affine Weyl group associated to the spherical Weyl group of type $X_n$, that is a semidirect product of the Weyl group of type $X_n$, called $W_{X_n}$, with a discrete group of translations. Starting with a compact or complex simple Lie group of type $X_n$ an explicit description of the Kac-Moody algebras of type $\widetilde{X}_n$ can be constructed via the intermediate step of the loop algebra, of $X_n$, that is the infinite dimensional Lie algebra consisting of maps 
$f:S^1\longrightarrow \mathfrak{g}$, satisfying some regularity conditions.

\noindent Then the affine Kac-Moody algebras of type $\widetilde{X}_n$ over the field $\mathbb{F}\in\{\mathbb{R}, \mathbb{C}\}$ are defined by 
\begin{displaymath}
\widehat{L}(\mathfrak{g}):=L(\mathfrak{g})\oplus c\mathbb{F}\oplus d\mathbb{F}
\end{displaymath}

\noindent such that $c$ is a central element and $d$ a derivation~\cite{Kac90,Terng95,Carter05,HPTT,Heintze09,Freyn09}.

Because of this close relationship between the simple Lie algebras of type $X$ and the associated affine Kac-Moody algebras of type $\widetilde{X}$, the linear representation of the twin city of type $\widetilde{X}$ relies in each case on the description of the spherical building of type $X$. 

As detailed in~\cite{Garrett97}, for the spherical buildings there exist the following descriptions:

\begin{enumerate}
 \item Buildings of type $A_n$ correspond to flag complexes of arbitrary flags.
 \item Buildings of type $C_n$ and $B_n$ correspond to flag complexes of subspaces isotropic with respect to some invariant form.
 \item Buildings of type $D_n$ correspond to the oriflamme complex, a slightly more complicated nesting construction of subspaces isotropic with respect to some invariant form. 
\end{enumerate}

To the knowledge of the author there are no linear representations for the buildings of the exceptional types $E_{n, n=6,7,8}$, $F_4$ and $G_2$ in the literature.

To construct linear representations of the cities, we use the embedding of loop groups into the restricted general linear group of a polarized Hilbert space $H=H^+ \oplus H^-$~\cite{PressleySegal86}. Those groups act on a special type of Grassmannians, introduced by Mikio Sato in the context of integrable systems. They come in pairs corresponding to the two mutually orthogonal subspaces $H^+$ and $H^-$. They consists of subspaces $W\subset H$, that are sufficiently close to $H^+$ resp.\ $H^-$. Cities are constructed as simplicial complexes associated to the space of flags satisfying a certain periodicity condition. In the case of groups of type $\widetilde{A}_n$, we get full flags. In the case of the groups of type $\widetilde{C}_n$ we need only special isotropic subspaces. In the case of types $\widetilde{B}_n$ and $\widetilde{D}_n$ we also use classes of isotropic subspaces, but apply a more complicated nesting procedure, corresponding roughly to the oriflamme construction used for finite dimensional buildings for groups of type $D_n$. 

We describe the content of the article in more details:
In section \ref{Grassmanians_and_periodic_flags} we review the functional analytic basics of our construction and describe the Sato Grassmannians and the periodic flag varieties.
The following three sections are each devoted to one class of groups, ordered in ascending complexity: unitary groups $\widetilde{A}_n$, symplectic groups $\widetilde{C}_n$ and orthogonal groups $\widetilde{B}_n$ and $\widetilde{D}_n$. Each of those sections consists of three parts. The first reviews the linear representations of the finite dimensional spherical buildings of type $X_n$, the second one is devoted to the construction of the affine city of type $\widetilde{X}_n$ and the third one studies the twinning.

For all other Kac-Moody groups the explicit construction of linear representations is still missing. As there are abstract, linear representations on Hilbert spaces for the nontwisted Kac-Moody groups of exceptional types~\cite{PressleySegal86}, it is clear that one can make this explicit. Hence a construction of the linear representation should be within reach. Nevertheless it is not evident that this description allows for a good geometric interpretation.  For the twisted affine Kac-Moody algebras there are explicit representations known in products of the eigenspaces of the eigenvalues of the diagram automorphisms. This picture extends to the associated Kac-Moody groups. Hence we conjecture also in this case the existence of linear representations of the associated cities.

\section{Functional analytic basics:\\ Grassmannians and periodic flag varieties}
\label{Grassmanians_and_periodic_flags}

In this section our aim is to investigate the functional analytic characteristics of the set of subspaces  $\mathcal{SUB}(V)$ for an infinite dimensional vector space $V$ that is appropriate for the construction of affine twin cities. As the finite dimensional blueprint suggests, twin cities of type $\widetilde{A}_n$ will be constructed using the whole of $\mathcal{SUB}(V)$, twin cities of type $\widetilde{C}_n$ will be constructed as the flag complex of a subset of ``isotropic'' subspaces in $\mathcal{SUB}(V)$, for twin cities of types $\widetilde{B}_n$,$\widetilde{D}_n$ we will have to implement versions of the oriflamme construction, starting with ``isotropic subspaces''. Remark that the isotropy condition is a purely algebraic one.

The functional analytic basics split into two main themes:

\begin{enumerate}
 \item regularity issues to adapt the Grassmannians to affine Kac-Moody groups $\widehat{L}(G)$ of various regularity conditions,
 \item algebraic issues to encode basic algebraic properties common to all Kac-Moody groups.
\end{enumerate}

Concerning the algebraic issues, there are three crucial aspects:
\begin{enumerate}
 \item the minimal dimension of matrix representations over function rings,
 \item the structure of the functions rings, especially multiplication by $z$,
 \item the involution $z\mapsto \frac{1}{z}$.
\end{enumerate}

We will define our subspaces by comparison with two distinguished isometric orthogonal subspaces $H^+$ and $H^-$. Regularity issues are dealt with by prescribing this comparison. Swapping the roles of $H^+$ resp.\ $H^-$ we get two symmetric objects, thus implementing the involution $z\mapsto \frac{1}{z}$. Accordingly we define two Grassmannians $Gr^{+}(H)$ and $Gr^-(H)$.    All subspaces in those two Grassmannians are infinite dimensional. Hence we have to be careful to find a new notion, replacing the finite dimensional concept of dimension: this is the so\ndash called virtual dimension. 
This type of infinite dimensional Grassmannians was introduced by Mikio Sato~\cite{Sato83} to describe integrable systems; their theory is treated in detail in chapters 6, 7 and 8 of the monograph~\cite{PressleySegal86}.  Compare also~\cite{SegalWilson85} and~\cite{Freyn09}.   We review the foundations and refer to~\cite{PressleySegal86} for further details.
The structure of the function ring will be implemented by the introduction of a shift operator. A wandering subspace for the shift operator connects to the matrix representation; its dimension reflects the minimal dimension of matrix representations.

We start by reviewing the representations of Kac-Moody groups as extensions of loop groups. Let $\mathfrak{g}$ be a simple complex or compact Lie algebra. Define the associated loop algebra by 
\begin{displaymath}
L(\mathfrak{g}):=\{f: \mathbb{R}\longrightarrow \mathfrak{g}| f(t+2\pi)=f(t),\, f \textrm{ satisfies some regularity condition}\}\, .
\end{displaymath}

\noindent The appended Kac-Moody algebra is defined by 
\begin{displaymath}
\widehat{L}(\mathfrak{g})=L(\mathfrak{g})\oplus \mathbb{F}c\oplus \mathbb{F}d\, ,
 \end{displaymath}
 where $d$ acts on $L(\mathfrak{g})$ as a derivative and $c$ is a central element. Hence $[d,f]=f'$, $[c,d]=[c,f]=0$ and $[f,g]=[f,g]_0+\omega(f,g)c$, where $[f,g]_0$ denotes the Lie bracket of $L(\mathfrak{g})$ and $\omega$ is a certain antisymmetric $2$\ndash form. 

 The minimal realization is the algebraic Kac-Moody algebra, corresponding to a Lie algebra of polynomial maps. All further Kac-Moody algebras arise as completions: on the one hand, understanding this algebra as an extension of a Lie algebra over the abstract polynomial ring $\mathbb{F}[t,t^{-1}]$, we can turn to the formal completion and study the resulting Kac-Moody algebras. On the other hand, taking the point of view of polynomial maps on $S^1$ resp.\ $\mathbb{C}^*$, we get a great variety of ``analytic'' completions and Kac-Moody algebras associated to them~\cite{Tits84}.

Correspondingly there is a variety of Kac-Moody groups. They are constructed as torus extensions $\widehat{L}(G)$ of loop groups 

\begin{displaymath}
L(G)=\left\{f:S^1\longrightarrow G\mid f \textrm{ satisfies some regularity condition}\right\},
 \end{displaymath}

for a simple complex or compact Lie group $G$. The $2$\ndash dimensional extension consists of a central $S^1$\ndash\ respective $\mathbb{C}^*$\ndash extension, corresponding to the $c$\ndash term of the Kac-Moody algebra and then a semidirect product with $S^1$ resp.\ $\mathbb{C}^*$ whose action on $L(G)$ is defined to be a shift of the argument: $w\cdot f(t)=f(t+w)$. For this shift to be well defined we need in the complex setting the elements $f\in L(G)$ to be defined on the whole of $\mathbb{C}^*$. This implies strong restrictions on the possible regularity conditions. For details see~\cite{PressleySegal86,Popescu05,Freyn09,Khesin09} and various other sources.

Some common regularity conditions are
continuous loops $LG$, $k$\ndash differentiable loops $L^kG$, smooth loops $L^{\infty}G$, real analytic or complex analytic loops $MG$ resp.\ $A_nG$ on $X\in \{\mathbb{C}^*, A_n:=\{z\in \mathbb{C}|e^{-n}\leq |z|\leq e^n\}\}$ (for the last two cases to make sense, we need $G$ to be a complex Lie group).
Furthermore there is the group of algebraic loops
$$L_{\textrm{alg}}G:=\{f \in L(G)|\textrm{ has a finite Fourier expansion}\},$$
where the Fourier expansion is defined via the adjoint representation of $G$. If $G$ is complex, we can identify this group with a group of matrix\ndash valued Laurent polynomials~\cite{PressleySegal86}. By construction, this group is isomorphic to the group $G(\mathbb{C}[t,t^{-1}])$, the realization of the affine algebraic group scheme corresponding to the Lie group $G$ over the ring $\mathbb{C}[t,t^{-1}]$ (for the definition see~\cite{AbramenkoBrown08} or~\cite{Waterhouse79}). Remark that $G(\mathbb{C}[t,t^{-1}])$ is the group acting in a natural way on a twin building~\cite{Ronan03,Remy02}.
If $G$ is semisimple compact or complex, then $L_{\textrm{alg}}G$ is dense in the group of continuous loops $LG$~\cite{PressleySegal86}, chapter~3.5. Hence $\widehat{L_{\textrm{alg}}G}$ is dense in $\widehat{L}(G)$. Thus the Kac-Moody groups $\widehat{L}(G)$ are completions
of the corresponding algebraic Kac-Moody groups $\widehat{L_{alg}G}$.

Loop groups embed naturally in the general linear group of Hilbert spaces:
Let $H^{n}= L^2(S^1, \mathbb{C}^n)$ denote the separable Hilbert space of square summable functions on $S^1$ with values in $\mathbb{C}^n$. 
Suppose a Lie group $G$ embeds in the general linear group $GL(n, \mathbb{C})$ acting on $\mathbb{C}^n$. Then a group $L(G)$ embeds into $L(GL(n))$, which is itself a subgroup of $GL(H^n)$~\cite{PressleySegal86}. 

 Let $H= H^{++} \oplus H^0 \oplus H^{--}$ be a polarization, that is a decomposition into the eigenspaces with positive eigenvalues $H^{++}$, zero eigenvalues $H^0$ and negative eigenvalues $H^{--}$ of some operator acting on $H$. We well use the one induced by the operator of $-i\frac{d}{d\theta}$. In the Fourier  basis, functions in $H^n$ are represented by a power series in $e^{i\theta}$:
\begin{displaymath}
 f(\theta)=\sum_n a_n e^{i \theta}
\end{displaymath} 
In this basis $H^{++}$ consists of all functions, such that the all notrivial terms in the Fourier series are positiv , $H^0$ consists of constant functions  and $H^{--}$ consists of all function, such that the only nontrivial terms in the Fourier series are negative. Each eigenspace has dimension $n$.  Set $H^+ = H^{++} \oplus H^0$ and $H^-= H^0 \oplus H^{--}$. $H^+$ consists of functions, that are boundary functions for holomorphic functions on the unit disc. If $n=1$, this space is the Hardy space $H^2$.

While a Hilbert space --- as any finite dimensional vector space --- has no distinguished subspaces, the additional structure of a polarization breaks this homogeneity. The second important structure element is the shift operator:

\begin{definition}[shift operator]
A surjective isometry $G: H^n\longrightarrow H^n$ is called a positive shift operator iff $G(H^{+})= H^{++}$. It is called a negative shift operator iff $G(H^{-})= H^{--}$.
\end{definition}

We know furthermore that $H^0$ is an orthogonal complement to $H^{++}$ in $H^+$. Hence the lemma of Wold-Kolmogorov~\cite{Nikolski02a}, lemma 1.5.1.\ tells us that $H^0$ is a wandering subspace for $G$ and we get the decomposition 

$$H^+=\underbrace{\left( \sum_{k\geq 0}\oplus G^k H^{0} \right)}_{=:H^+_{\textrm{fin}}}\oplus \underbrace{\left(\bigcap_{k\geq 0} G^k H^0\right)}_{=:H^{+}_{\textrm{inf}}} =H^+_{\textrm{fin}}\oplus H^{+}_{\textrm{inf}}\, .$$ 

 Suppose $x\in H^+_{\textrm{fin}}$, then there is a minimal $k\geq 0$ such that $x \in \sum_{i\leq k} G^i(H_0)$. $k$ is called the degree of $x$. Thus $H^+_{\textrm{fin}}$ corresponds to the subspace of polynomials in $e^{i\theta}$ of $H^+$ while $H_{\textrm{inf}}$ is the complement. 
Surjectivity of $G$ yields that $G(H^{--})=H^-$. Hence $G^{-1}:H^n\longrightarrow H^n$ satisfies $G(H^{-})= H^{--}$. Hence $G^{-1}$ is a negative shift operator. The shift operator maps the eigenspace with eigenvalue $n$ onto the eigenspace with eigenvalue $n+1$. 

\begin{example} 
Using the realization $H = L^2(S^1, \mathbb{C}^n)$ we can define the shift operator to be the multiplication operator $M_z f= zf$ 
($z=e^{i\theta}$ on $S^1$).
\end{example}

In Fourier Analysis, this operator is referred to as shift operator.  $H^0$ is called a wandering subspace for the grading operator~\cite{Nikolski02a,Nikolski02b}.
Let us note the following theorem~\cite{PressleySegal86}, theorem 6.1.1:

\begin{theorem}
 The commutant of $M_z$ in $GL(H^{(n)})$ is the group $L_{meas}GL_n(\mathbb{C})$ of bounded measurable maps $S^1 \longrightarrow GL_n(\mathbb{C})$.
\end{theorem}

\noindent Following~\cite{PressleySegal86}, definition 7.1, the positive Grassmannian is defined as follows:

\begin{definition}[positive Grassmannian]
\label{positivegrassmanian}
The positive Grassmannian $Gr^{+}(H)$ is the set of all closed subspaces $W$ of $H$ such that
\begin{enumerate}
	\item the orthogonal projection $pr_+: W \longrightarrow H^+$ is a Fredholm operator,
	\item the orthogonal projection $pr_{--}: W \longrightarrow H^{--}$ is a compact operator.
\end{enumerate}
\end{definition}

\begin{remark}
Sometimes (see for example~\cite{PressleySegal86}) the second condition is replaced by the requirement for $pr_{--}: W \longrightarrow H^{--}$ to be a Hilbert-Schmidt operator. Hilbert-Schmidt operators form an ideal in the set of compact operators. The square of the ideal of Hilbert-Schmidt operators consists of trace class operators. Restriction of the theory to those Hilbert-Schmidt Grassmannians is necessary if one is interested in the construction of certain determinant bundles. 
\end{remark}

\noindent One can assign to each element of $Gr^+(H)$ an integer called the virtual dimension:

\begin{definition}[virtual dimension]
\label{virtualdimension}
Let $W\in Gr^+(H)$.  The  virtual dimension of $W$ is defined by
\begin{displaymath}
\nu (W) =\textrm{dim}(\textrm{ker } pr_+) -\textrm{dim}(\textrm{coker } pr_{+})\, .
\end{displaymath}
\end{definition}

In many situations, the condition on $pr_{--}$ to be a compact operator is too weak and the resulting Grassmannian too big. Many conditions on the operator $pr_{--}$ correspond to regularity conditions of functions - geometrically, they describe how ``far'' away a subspace $W$ can be maximally from $H^+$ (resp.\ $H^-$). This leads to various versions of restricted Grassmannians.  They are defined in a way that loop groups (resp.\ Kac-Moody groups) satisfying usual regularity conditions act nicely on them.  The smallest Grassmannian is the algebraic one:

\begin{definition}[positive algebraic Grassmannian]
\label{algebraicgrassmanian}    
The positive algebraic Grassmannian $Gr_0^+(H)\subset Gr^+(H)$ consists of subspaces $W\in Gr^+(H)$ such that 
\begin{displaymath}
z^k H^+ \subset W \subset z^{-k} H^+\, .
 \end{displaymath}
 \end{definition}

Using the explicit description $ H=L^2(S^1, \mathbb{C}^n)$, $Gr_0^+(H)$ consists exactly of the elements $W\in Gr(H)$ such that the images of $pr_{--}:W \longrightarrow H^{--}$ and $pr_+: W^{\perp} \longrightarrow H^+$ are polynomials. Stated in an operator theoretic way it consists of those operators, such that there is a $k>0$ such that for $|i-j|>kn$ the coefficient $a_{ij}=0$. In the language of $n$\ndash Laurent operators this is equivalent to $W$ being defined by an algebraic loop of degree less than $k$~\cite{PressleySegal86} and~\cite{Gohberg03}.

$Gr_0^+(H)$ is the union of finite dimensional Grassmannians.

\noindent For the definition of similar subspaces -- the rational Grassmannian $Gr_1^+(H)$, $Gr_{\omega}^+(H)$, the Hilbert-Schmidt Grassmannian $GR_{HS}^+(H)$ and the smooth Grassmannian $Gr_{\infty}^+(H)$ see~\cite{PressleySegal86}. For the applications in Kac-Moody geometry~\cite{Freyn10d}  we need two new regularity conditions: To describe the building in the situation of polar actions on Hilbert spaces as considered by Chuu-Lian Terng~\cite{Terng95}, we use Sobolev Grassmannians:

\begin{definition}[positive $H^1$\ndash Grassmannian]

The $H^1$\ndash positive Sobolev Grassmannian $Gr_{H_1}^+(H)$ consists of the graphs of operators $T: H_s \longrightarrow H_s^{\perp}$ whose entries $T_{pq}$ satisfy: $|T_{pq}|{(p-q)^2} < \infty $.   
\end{definition}

\noindent Furthermore in connection with Kac-Moody symmetric spaces that are tame Fr\'echet manifolds~\cite{Freyn10a,Freyn09}, we need a tame Fr\'echet Grassmannian:

\begin{definition}[positive tame Fr\'echet Grassmannian]
\label{tamefrechetgrassmanian}
The positive tame Fr\'echet Grassmannian $Gr_{t}^+(H)$ consists of the graphs of operators $T: H_s \longrightarrow H_s^{\perp}$ whose entries $T_{pq}$ are exponentially decreasing: $|T_{pq}|e^{(p-q)n} < \infty \hspace{3pt}\forall n$.
\end{definition}

Compare this definition with the notion of exponential weights developed in~\cite{GoodmanWallach84}.

\noindent In the language of $n$\ndash Laurent operators this condition is equivalent to the loop being in $MG$, i.e.\ the restriction to holomorphic loops defined on the whole of $\mathbb{C}^*$~\cite{Gohberg03}.

\begin{definition}[reduced Grassmannian]
\label{reducedgrassmanian}
The reduced positive Grassmannian of regularity $X$ denoted $Gr^{n,+}_X(H^n)$ consists of subspaces $W\subset Gr^+(H^n)$ that are invariant with respect to the Grading operator (i.e. such that $G(W) \subset W$ or explicitly $zW\subset W$).
\end{definition}

The definition of the other types of reduced Grassmannians, especially reduced algebraic Grassmannians, reduced $H^1$- and reduced tame Grassmannians is self explaining. 

The following theorem~\cite{PressleySegal86}, theorem 8.3.2, shows this to be the correct notion to work well with the action of loop groups.
Following ~\cite{PressleySegal86} we define $GR_{HS}^+(H)$ to consist of all subspaces in $Gr^+(H)$ such that the orthogonal projection $pr_{--}: W \longrightarrow H^{--}$ is a Hilbert-Schmidt operator.

\begin{theorem}
The group of $\frac{1}{2}$\ndash differentiable Sobolev loops $L_{\frac{1}{2}}U_n$ acts transitively on the Hilbert-Schmidt Grassmannian $Gr_{HS}^{n,+}(H^n)$ and the isotropy group of $H^+$ is the group $U_n$ of constant loops. 
\end{theorem}

This theorem yields the equivalences $\Omega_{\frac{1}{2}}U_n = Gr^{n}(H)$ and $\Omega_{alg}U_n= Gr_0^{n}$. Similar statements hold for $Gr_1^{n,+}(H)$, $Gr_{\omega}^{n,+}(H)$, $Gr_{\infty}^{n,+}(H)$  and $Gr_{t}^{n,+}(H)$; 
for $Gr_1^{n,+}(H)$, $Gr_{\omega}^{n,+}(H)$, $Gr_{\infty}^{n,+}(H)$ a proof can be found in~\cite{PressleySegal86}; this proof adapts to the case of $Gr_{t}^{n,+}(H)$ straight forwardly.

Let $(e_1, \dots, e_n)$ be a $\mathbb{C}$\ndash basis of $H^0$ and  define  $R=L^2(S^1, \mathbb{C})$. Then $(e_1, \dots, e_n)$ is an $R$\ndash basis of $H^n$. Let  $f\in L_{X}U_n$ where $X$ denotes one of the introduced regularity conditions. Then --- as $L_{X}U_n$ is a group, $f$ has an inverse element ---  $\left(f(e_1), \dots, f(e_n)\right)$ is also an $R$\ndash basis of $H^n$. Define $R^+:=\{f\in R| \textrm{ the negative Fourier components of } f \textrm{ vanish}\}$. Then $(e_1, \dots, e_n)$ is an $R^+$\ndash basis of $H^+$ and the $R^+$\ndash span of $\left(f(e_1), \dots, f(e_n)\right)$ is an element in $Gr_X^{n,+}(H^n)$, the periodic Grassmannian of regularity $X$.

\noindent The next step is the definition of the flag varieties:
Following again~\cite{PressleySegal86}, definition 8.7.5, we define (full) periodic flags:

\begin{definition}[full positive periodic flag manifold]

The full positive periodic flag manifold $Fl^{n,+}$ consists of all sequences $W_k, k \in \mathbb{Z}$ of subspaces in $H^{n}$ such that

\begin{enumerate}
	\item $W_k \subset Gr^{n,+}H$,
	\item $W_{k+1} \subset W_k \forall k$ and $\textrm{dim}(W_k/W_{k+1})=1$, 
	\item $z W_k = W_{k+n}$.
\end{enumerate}
\end{definition}

As elements of a flag satisfy $z W_k = W_{k+n}\subset W_k$, all elements of a flag are taken from $Gr^{n,+}$.

Let $\{e_1,\dots, e_n\}$ be a basis of $V^n\simeq \mathbb{C}^n$ and $V_i:=\textrm{span}\langle e_{i+1}, \dots, e_n\rangle$. 
We define the positive normal flag to be the flag $\{W_{k'}\}_{k'\in \mathbb{Z}}$ such that $W_{k'}=W_{kn+l}:=z^k W_{l}$  for $k'=kn+l, k\in \mathbb{Z}, l\in \{0, \dots, n-1\}$ and $W_l:=\{f:\mathbb{C}\longrightarrow V^n|f \textrm{ is holomorphic and } f(0)\subset V_i\}$.

To define the manifolds of partial periodic flags the virtual dimension is used. It behaves well with respect to inclusions:

As the index is additive under composition of Fredholm operators, i.e.\ for two Fredholm operators $A$ and $B$ we have $\textrm{ind}(AB)= \textrm{ind}(A) + \textrm{ind}(B)$ (see~\cite{Conway90}, chapter 11), one calculates for the virtual dimension: 
Let $W \subset W' \in Gr^+(n)$ and $\textrm{dim}\left( W/W'\right)=l$. Then $\nu(W')=\nu(W)+l$. For $W\in \textrm{Gr}^{n,+}$ an important special case is $\nu(W_{k+n})=\nu(zW_{k})= \nu(W_{k})-n$.

We now aim for the definition of partial flags.
For this purpose let $\mathcal{K}\subset \mathbb{Z}$ such that with $k\in \mathcal{K}$ also $k+nl\in \mathcal{K}, \forall l\in \mathbb{Z}$. Furthermore set  $\mathcal{K}_n:= \mathcal{K}\cap \{1,\dots, n\}$ and define $m_{\mathcal{K}}:=\#\mathcal{K}_n$.

\begin{definition}[partial positive periodic flag manifold]

The positive periodic flag manifold $Fl_{\mathcal{K}}^{n,+}$ consists of all flags $\{W_k\}, k \in \mathbb{Z}$  in $H^{n}$ such that
\begin{enumerate}
	\item $W_k \subset \textrm{Gr}^+H$,
	\item $W_{k+1} \subset W_k \forall k$, 
	\item $zW_k = W_{k+n}$.
	\item For every flag $\{W_k\}_{k\in \mathbb{Z}}$ the map $\nu:(\{W_k\}_{k\in \mathbb{Z}})\longrightarrow \mathbb{Z}$ mapping every subspace $W_k$ to its virtual dimension is a surjection onto $\mathcal{K}$.
\end{enumerate}
\end{definition}

This definition contains the one for full flags by using $\mathcal{K}=\mathbb{Z}$. In contrast if $m_{\mathcal{K}}=1$, we call a flag elementary. Elementary flags are in bijection with elements of $Gr^{n,+}(H)$ using the identification $Gr^{n,+}(H)\ni W_0 \leftrightarrow \{z^k W_0\}_{k\in \mathbb{Z}}\in Fl_{\nu(W_0)+n\mathbb{Z}}^{n,+}$.

\noindent We normalize our notation by requiring that $W_0$ is the element with the lowest nonnegative virtual dimension. In the case of full flags this is $0$, in the case of partial flags it may be between $0$ and $n\hspace{-1pt}-\hspace{-1pt}1$.

\noindent The description of the structure of $Fl^n$ shows a direct connection with results of~\cite{Mitchell88}. The crucial point is that all manifolds of periodic flags fiber over the Grassmannian with the fibers being finite (partial) flag manifolds.

For the maximal flag manifold this fiber structure can be described  by the exact sequence:
$$ 1 \longrightarrow Fl(\mathbb C^n) \longrightarrow Fl^{(n)} \longrightarrow Gr^{(n)} \longrightarrow 1\,.$$ 

\noindent To see this, note that the periodicity condition determines a whole maximal flag as soon as the spaces $\{W_0, W_1, \dots, W_{n-1}\}$ are specified. After the identification of $W_0/W_{n}\simeq \mathbb{C}^n$ the spaces $\{W_1, \dots, W_{n-1}\}$ determine a maximal flag in $\mathbb{C}^n$.

To construct a similar fibration for manifolds of partial flags the fiber has to be replaced by manifolds of partial flags. In doing so we have to construct the fibers in a way that the resulting flags get the correct virtual dimensions:

So let $\{W\}_{k\in \mathbb{Z}}$ be a flag in $Fl_{\mathcal{K}}^{n,+}$. Take the projection:  $\{W_k\}_{k\in \mathbb{Z}} \longrightarrow W_0$. Again $W_0/zW_0\simeq \mathbb{C}^n$. The intersections $U_l:=W_l\cap W_0/zW_0,l\in \{1,m_{\mathcal{K}}\}$ determine the flag completely. They describe a finite dimensional partial flag $U_{m_{\mathcal{K}}-1} \subset \dots \subset {U_1}$  in $zW_0/W_0\simeq \mathbb{C}^n$. The dimensions are $\textrm{dim}(U_l)=\textrm{virt dim}(W_l)-\textrm{virt dim}(zW_0)$.

\noindent This shows that $Fl^{n,+}_{\mathcal{K}}$ fibers over $Gr^{n,+}$, the fiber being $Fl_{\kappa_1, \dots,\kappa_{m_{\mathcal{K}}}}(\mathbb C^n)$ where $\kappa_i:=k_i+\textrm{virt dim}(W_0)$.

$$ 1 \longrightarrow Fl_{\kappa_1, \dots,\kappa_{m_{\mathcal{K}}}}(\mathbb C^n) \longrightarrow Fl_{\mathcal{K}}^{(n)} \longrightarrow Gr^{(n)} \longrightarrow 1\,.$$

Up to now we focused our attention to the positive part.
A completely symmetric theory can be developed for $H^-$:

\begin{definition}[Negative Grassmannian]

The negative Grassmannian $Gr^{-}(H)$ is the set of all closed subspaces $W$ of $H$ such that
\begin{enumerate}
   \item the orthogonal projection $pr_-: W\longrightarrow H^-$ is a Fredholm operator,
   \item the orthogonal projection $pr_{++}: W \longrightarrow H^{++}$ is a compact operator.
\end{enumerate}

\end{definition}

The subspaces $Gr_0^-(H)$, $Gr_1^-(H)$, $Gr_{\omega}^-(H)$, $Gr_{\infty}^-(H)$, $Gr_{t}^-(H)$, $Gr_0^{n,-}(H)$, $Gr_1^{n,-}(H)$, $Gr_{\omega}^{n,-}(H)$, $Gr_{\infty}^{n,-}(H)$ and $Gr_{t}^{n,-}(H)$ are defined as in the positive case.

We define full negative flags by:

\begin{definition}[full negative periodic flag manifold]

The full negative periodic flag manifold $Fl^{-,n}$ consists of all sequences $V_k, k \in \mathbb{Z}$ of subspaces in $H^n$ such that
\begin{enumerate}
	\item $W_k \subset \textrm{Gr}^{-}(H)$,
	\item $W_{k+1} \subset W_{k} \forall k$, $\textrm{dim}(W_k/W_{k+1})=1$,
	\item $\frac{1}{z} W_k= W_{k+n}$.
\end{enumerate}
\end{definition}

The negative normal flag and partial flags are defined similarly as in the case of positive flags.

\section{The special linear groups: Type \protect\boldmath$\widetilde{A}_n$}
\label{The_special_linear_groups_Type_A_n}

\subsection{The finite dimensional blueprint: type $A_n$}

The complex simple Lie group of type $A_{n-1}$ is $SL(n,\mathbb{C})$, its compact real form is $SU(n)$, and its split real form is $SL(n,\mathbb{R})$.
$SL(n; \mathbb{K})$ for $\mathbb{K}\in \{\mathbb{R}, \mathbb{C}\}$ acts on $V_{\mathbb{K}}^n\cong \mathbb{K}^n$.  

Its Dynkin diagram is

\centerline{
\begin{picture}(150,20)
\multiput(36,9)(30,0){4}{\circle{4}} 
 \put(33,0){\small $\alpha_1$\normalsize}\put(63,0){\small $\alpha_2$\normalsize}\put(91,0){\small $\alpha_{n-1}$\normalsize}\put(121,0){\small $\alpha_{n}$\normalsize} 
\multiput(42,9)(60,0){2}{\line(22,0){16}}
\put(74.5,8){$\dots$}
\end{picture} 
}

The Weyl group is the symmetric group $Sym(n)$. Hence it acts naturally by permutation on any basis $\{e_1, \dots e_n\}$ of $V_{\mathbb{K}}^n$.

We define 
\begin{displaymath}
\mathcal{F}_{\mathbb{K}}=\{F=(V_1\subset V_2 \subset \dots \subset V_i)|V_l\subset V_{\mathbb{K}}\, \quad\textrm{for }l=1,\dots ,i\}
\end{displaymath}
to be the set of all flags in $V^n_{\mathbb{K}}$. 

 A flag $F\subset \mathcal{F}_{\mathbb{K}}$ is said to be of type $\mathcal{K}$ for $\mathcal{K}=(k_1, \dots, k_j)\subset (1, \dots, n)$ iff $\textrm{dim}(V_i) = k_i, i=1, \dots, j$. We denote by $\mathcal{F}_{\mathcal{K}}$ the set of flags of type $\mathcal{K}$. Hence 

\begin{displaymath}
\mathcal{F}_{\mathbb{K}}=\bigcup_{\mathcal{K}\in \mathbb{P}(n)} \mathcal{F}_{\mathbb{K},\mathcal{K}}\, .
\end{displaymath}

 A flag is called maximal if $\mathcal{K}=\{1, \dots, n\}$ or equivalently $\textrm{dim}_{\mathbb{K}}V_i/V_{i+1}=1$ for $i\in \{1, \dots, n-1\}$ otherwise it is called partial.

The $SL(n)$\ndash action on $V$ induces an action on $\mathcal{F}$ in a natural way:

$$SL(n,\mathbb{K}):\,\mathcal{F}\,\longrightarrow\, \mathcal{F} \qquad
g\cdot(V_1\subset \dots \subset V_i)\mapsto (g\cdot V_1\subset \dots \subset g \cdot V_i)\, .
$$

Let $F=(V_1\subset V_2 \subset \dots \subset V_j)$ be a flag. The stabilizer subgroup $P_F\subset SL(n,\mathbb{C})$ of $F$ consists of all $g\in SL(n, \mathbb{C})$ such that $g\cdot V_i=V_i$ for all $i$. 

It is elementary to check the following well-known facts:
\begin{enumerate}
\item The Borel subgroups of $SL(n, \mathbb{K})$ are exactly the stabilizer subgroups of maximal flags. The standard Borel subgroup $B_0$, i.e. the group of upper triangular matrices, stabilizes the standard flag 

\begin{displaymath}
F=(V_i=\textrm{span}(e_1, \dots, e_i) \textrm{ for } i=1, \dots, n)\, .
\end{displaymath}

\item The parabolic subgroups are stabilizers of (partial) flags.
\item The $SL(n, \mathbb{K})$\ndash action preserves the $\mathcal{K}$\ndash type of flags.
\item The $SL(n, \mathbb{K})$\ndash action is transitive on $\mathcal{F}_{\mathcal{K}}$.
\end{enumerate}

Hence there is a $SL(n,\mathbb{K})$\ndash equivariant bijection between parabolic subgroups of $SL(n,\mathbb{K})$ and flags in $V^n\simeq\mathbb{K}^n$ (remember that Borel subgroups are parabolic).

The building $\mathfrak{B}$ associated to $Sl(n, \mathbb{K})$ is isomorphic to the flag complex over the set of parabolic subgroups of $SL(n, {\mathbb{K}})$, with the order relation defined by reverse inclusion.
\begin{displaymath}
 P_1< P_2\qquad \Leftrightarrow\qquad P_1 \supset P_2\, .
\end{displaymath}
 
Identifying parabolic subgroups with the (partial) flags stabilized by them we can reformulate this in terms of flags:

We define a partial ordering on the set of flags by 
\begin{displaymath}
F_1 < F_2\qquad \Leftrightarrow \qquad F_1 \subset F_2
\end{displaymath}
as sets of subspaces. Geometrically we identify a flag $F=\{V_i\}_{i\in I}$ with a simplex $\Delta_F$ of dimension $|I|-1$. $\Delta_{F_1}$ is in the boundary of $\Delta_{F_2}$ iff $F_1 < F_2$. This simplicial complex associated to the set of all flags is called the flag complex. We have the following result:

\begin{lemma}
\label{spherical_building_type_A_n}
The flag complex is  isomorphic to the building. 
\end{lemma}

\begin{proof}For details see~\cite{Garrett97}.
\end{proof}

In the flag complex description apartments are characterized by frames:
 We define a frame $A_f$ for $V^n$ to be a set of $n$ $1$\ndash$\mathbb{K}$\ndash dimensional subspaces of $V^n$, denoted $U_i$, such that $V=\textrm{span} \langle U_1, \dots, U_n \rangle$.  A subspace $V_i\subset V^n$ is called subjacent to the frame $A_f$ iff there are elements $\{U_{i_1}, \dots, U_{i_l}\}\in A_f$ such that $V_i=\textrm{span}\langle U_{i_1}, \dots, U_{i_l}\rangle$.

The apartment $A$ corresponding to the frame $A_f$ consists of all simplices corresponding to flags subjacent to $A_f$. If a flag $F_1$ is subjacent to $A_F$ so is every flag $F_2\subset F_1$. Hence $A$ is a simplicial complex.

The Weyl group of $SL(n, \mathbb{K})$ is the symmetric group in $n$ letters. It is realized as the permutation group of the elements of any frame. This action induces an action on the subspaces subjacent to the frame and hence on the flags constructed of those subspaces and thus on the apartment.



\subsection{The affine city} 

In this section we will prove the following result:

\begin{theorem}
 The simplicial complex associated to the ordered set of all positive and (resp.\ negative) (partial) periodic flags is an affine city.  
\end{theorem}

It is the flag complex associated to the set of all elementary flags where the incidence relation between two elementary periodic flags is defined as follows:

\begin{definition}[Incident relation]
Let $W$ and $W'$ be two subspaces in $Gr_X^{+,n}$. The elementary flags $\{W'_{kn}\}:=z^{k} W'_l$ and $\{W'_{km}\}:=z^{m} W'$ are incident iff there are $l_1$ and $l_2$ such that $W \subset z^{l_1}W' \subset z^{l_1+l_2}W$.
\end{definition}

In this section we prove that the flag complex associated to the positive periodic flags is a city.

To prove this claim, we have to
\begin{enumerate}
 \item show that this complex is a chamber complex,
 \item define apartments,
 \item check the axioms required for the system of apartments. 
\end{enumerate}

\begin{definition}[Boundary]
Let $\{W_k\}$ be a flag. Then the set of flags $\{W'_l\}$ such that for each $l$ there is a $k$ such that $W'_l = W_k$ is denoted by $\{\leq W_k\}$ and is called the boundary of $\{W_k\}$.
\end{definition}

\begin{lemma}
\label{chambercomplex}
The flag complex of all periodic flags is an $n$\ndash dimensional simplicial complex, such that every simplex is contained in a simplex of maximal dimension.
\end{lemma}

\begin{cor}
\label{chambercomplex_cor}
It is not a chamber complex as the definition of a chamber complex includes connectedness. Hence every connected component is a chamber complex.
\end{cor}

\begin{proof}[Proof of lemma \ref{chambercomplex}] The complex of all (partial) periodic flags is a chambers complex as for every (partial) periodic flag $\{W_k\}_{k\subset \mathcal{K}}$ of type $\mathcal{K}$ and for every $\mathcal{K}' \subset \mathcal{K}$ 
there is a partial flag $\{W'_k\}_{k\in \mathcal{K}'}$ of type $\mathcal{K'}$ such that $\{W'_k\}_{k\in \mathcal{K}'} \leq \{W_k\}_{k\in \mathcal{K}}$.

To prove that it is a chamber complex we have to check, that every partial periodic flag may be completed to a maximal periodic flag. Choose a set $\mathcal{K}\not= \mathbb{Z}$ and let $\{W_k\}$ be a partial flag of type $\mathcal{K}$. As $\mathcal{K}\not=\mathbb{Z}$ there is k such that $dim(W_k/W_{k+i})\geq 2$ for two consecutive elements $k, k+i\subset \mathcal{K}$. Hence we can find an element $W_{k+j}$ such that $W_{k+i}\subset W_{k+j}\subset W_k$. Extend $W_{k+j}$ to the elementary flag $z^l W_{k+j}, l\in \mathbb{Z}$. By periodicity the elements of this flag satisfy $W_{ln+k+i}\subset W_{ln+k+j}\subset W_{ln+k}$.
 Hence the new flag is of type $\mathcal{K}'$ such that $m_{\mathcal{K}'}=m_{\mathcal{K}}+1$. We reiterate this procedure till after $s$ steps we arrive at $\mathcal{K}^s =\mathbb{Z}$. Thus we have constructed a full flag and hence proved that each simplex is contained in the boundary of a chamber.
\end{proof}

Following the finite dimensional blueprint we describe apartments using (periodic) frames:

\begin{definition}[frame]
A (periodic) frame is a sequence of subspaces $\{U_k\}_{k\in \mathbb{Z}}\subset H^n$ such that $U_{k+n}=z U_k$ and $H^n=\textrm{span}_k U_k$.
\end{definition}

Let $\{W_k\}_{k\in \mathbb{Z}}$ denote a full flag. A frame such that $W_k:=W_{k+1}\oplus U_k$ will be called normal with respect to $\{W_k\}$. Clearly  any flag is normal to  frames.

\begin{definition}[admissible permutation]

We call a permutation $\pi: \{U_k\}_{k \in\mathbb{Z}}\longrightarrow \{U_k\}_{k \in\mathbb{Z}}$ admissible if $\pi(U_{k+n})= \pi(U_k)+n$.
\end{definition}

\begin{definition}[affine Weyl group associated to $\{U_k\}_{k \in\mathbb{Z}}$]
The affine Weyl group $W_{\textrm{aff}}$ is defined to be the group of admissible permutations of any frame $\{U_k\}_{k \in\mathbb{Z}}$.
\end{definition}

$W_{\textrm{aff}}$ is independent of the choice of the periodic frame $\{U_k\}_{k\in \mathbb Z}$. This definition coincides with the usual definition for the affine Weyl group of the Kac-Moody algebra of type $\widetilde{A}_n$~\cite{Kac90} and~\cite{PressleySegal86}. 
The Dynkin diagram of type $\widetilde{A}_n$ is

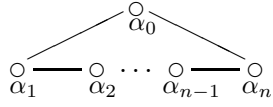
\centerline{
\begin{picture}(140,45)
\multiput(36,9)(30,0){4}{\circle{5}} 
 \put(33,0){\small $\alpha_1$\normalsize}\put(63,0){\small $\alpha_2$\normalsize}\put(91,0){\small $\alpha_{n-1}$\normalsize}\put(121,0){\small $\alpha_{n}$\normalsize}    
\multiput(42,9)(60,0){2}{\line(22,0){18}}
\put(74.5,8){$\dots$}
\put(81,32){\circle{5}}\put(78,24){\small $\alpha_{0}$\normalsize}\put(40,13){\line(2,1){36}}\put(86,31){\line(2,-1){36}}
\end{picture} 
}
$W_{\textrm{aff}}$ is generated by the set of transformations $S:=\langle \alpha_{i}|i=0,\dots, n-1\rangle$ such that 
\begin{alignat*}{1}
\alpha_{i}(U_{ln+i})&= U_{ln+i+1},\\
\alpha_{i}(U_{ln+i+1})&= U_{ln+i},\\
\alpha_{i}(U_{ln+k})&= U_{ln+k}, k\not= i, i+1.
\end{alignat*}

\noindent We have $W_{\textrm{aff}}\simeq \mathbb{Z}^{n-1} \rtimes \textrm{Sym}(n)$. Let us investigate the action of the two factors in the decomposition $W_{\textrm{aff}}\simeq \mathbb{Z}^{n-1} \rtimes \textrm{Sym}(n)$:

An element $\pi\in  \textrm{Sym}(n)$ acts on $\{U_{ln+k}\}$ via permutation of $k$: 
$$\pi (\{U_{ln+k}\})=\{U_{ln+\pi(k)}\}\,.$$

An element $A=(a_1, \dots, a_{n-1})\in \mathbb{Z}^{n-1}$ acts on $\{U_{ln+1}, U_{ln+2}, \dots, U_{ln+n-1} \}$ via
\small\begin{displaymath}
A(\{U_{ln+1}, U_{ln+2}, \dots, U_{ln+n-1} \})= \{U_{(l-a_1)n+1}, U_{(l+a_1-a_2)n+2}, \dots, U_{(l+a_{n-1})n+n-1} \}.
\end{displaymath}\normalsize

Admissible permutations induce maps on the flag manifolds. An admissible permutation $w: \{U_l\}\longrightarrow w(U_{l})$ of the frame $\{U_l\}$ maps the flag 
\begin{displaymath}\{W_k\}:= \{\textrm{span}_{l\geq k} U_l\}\end{displaymath}
onto the flag \begin{displaymath}w(\{W_k\})= \{\textrm{span}_{l\geq k} w(U_l)\}.\end{displaymath} 
For each element $W_k$ in a flag $\{W_k\}$ one can check: $\nu(W_k)= \nu (w(W_k))$.

\begin{lemma}
Let $\{U_l\}$ be any frame, represent $W$ as the group of admissible permutations of $U_n$.
The induced action of the Weyl group on the restricted Grassmannian preserves the virtual dimension.
Furthermore the induced action on the space of all flags $\mathcal{F}$ preserves the $\mathcal{K}$\ndash type of flags.
\end{lemma} 

In other words: The $W_{\textrm{aff}}$\ndash orbit of a flag $\{W_{k}\}\in Fl_{\mathcal{K}}^{n,+}$ is completely contained in $Fl_{\mathcal{K}}^{n,+}$.

\begin{proof}
The second part of this lemma is a direct consequence of the first part. For the
first part it is sufficient to show that it is satisfied for all generators
$\alpha_{i}$; to verify this, one has to distinguish between two cases: Let $\alpha_{i}$
interchange $U_{i+nl}$ and $U_{i+1+nl}, l\in \mathbb{Z}$. If
with  $U_{i+nl}\subset W_k$ also $U_{i+1+nl} \subset W_k$, then $W_k$ is stable
under $\alpha_{i}$, so also the virtual dimension is stable; if there is
$U_{i+1+nl} \subset W_k$ such that $U_{i+nl} \not\subset W_k$, then
$\textrm{dim}\left(W_k/ W_k \cap \alpha_{i}(W_k) \right)= 1$. So 
$\textrm{virt dim}(W_k) = 1+\textrm{virt dim}(W_k \cap \alpha_{i}(W_k))$. 
Symmetrically, $\textrm{dim}\left(\alpha_{i}(W_k)/ (W_k \cap \alpha_{i}(W_k)) \right)= 1$  
leads to $\textrm{virt dim}(\alpha_{i}(W_k)) = 1+\textrm{virt dim}(W_k \cap \alpha_{i}(W_k))$, 
which shows the equivalence. 
\end{proof}

We now need a description of regularity conditions for frames:

\begin{lemma}
\label{regularity conditions for frames}
 Let $\{U_k\}$ be a periodic frame and let $Gr_x^{n,+}$ be a Grassmannian of some prescribed regularity. The following are equivalent:
\begin{enumerate}[(i)]
 \item The space $W_0:= \bigoplus_{i=0}^{\infty} \{U_i\}$ is in $Gr_x^{n,+}$.
 \item All spaces $W_k:= \bigoplus_{i=k}^{\infty} \{U_i\}$ are in $Gr_x^{n,+}$.
 \item Let $\pi:\{U_k\}\longrightarrow \{U_k\}$ be an admissible permutation. All spaces $W_{\pi, k}:= \bigoplus_{i=k}^{\infty} \{U_{\pi(i)}\}$ are in $Gr_x^{n,+}$.
\end{enumerate}
\end{lemma}

\begin{proof}
 The implications $(iii)\Rightarrow (ii)$ and $(ii)\Rightarrow (i)$ are trivial.
The implication $(i)\Rightarrow (ii)$ follows as $W_k=U_k\oplus W_{k+1}$. Hence $W_k\in Gr_x^{n,+}$ iff $W_{k+1}\in Gr_x^{n,+}$. The implication $(ii)\Rightarrow (iii)$ follows, as $W_{\pi, k}\cap W_k$ has finite codimension in $W_{\pi,k}$ and in $W_k$. Hence $W_k \in Gr_x^{n,+}$  is equivalent to $W_{\pi, k} \in Gr_x^{n,+}$. 
\end{proof}

\begin{definition}
 A frame $\{U_k\}$ satisfying one and hence all equivalent conditions of lemma \ref{regularity conditions for frames} is called a frame of regularity $X$.
\end{definition}

\begin{definition}[apartment]

Let $\{W_k\}$ be a full flag, $\{\leq \hspace{-3pt}W_k\}$ the boundary of $\{W_k\}$ and $\{U_k\}_{k \in\mathbb{Z}}$ a frame. The apartment $\mathcal{A}(\{W_k\}, U_k)$ consists of all flags that are admissible permutations of flags in $\{\leq W_k\}$ with respect to the frame $\{ \{U_k\}_{k \in\mathbb{Z}}\}$.
\end{definition}

An apartment is an abstract simplicial complex. It can be identified with a simplicial chamber complex, whose chambers correspond to full flags. As $W_{\textrm{aff}}$ is a Coxeter group that acts transitively on the chambers of the apartment, it is in fact a Coxeter complex. 

The apartments containing a given flag $\{W_k\}$ are in bijection with frames $\{U_k\}$ associated to this flag. 

For a pair of two flags of the same $\mathcal{K}$\ndash type $\{W_k\}$ and $\{W'_{k'}\}$, there exists an apartment containing both of them iff there is a frame $\{U_k\}$ for $\{W_k\}$, such that for some $w\in  W_{\textrm{aff}}$ $w(\{U_k\})$ is a frame for $\{W'_{k'}\}$. 

In contrast to the finite dimensional and even the infinite dimensional algebraic situation, it is not true that every pair of full flags is contained in a joint apartment. Now we have:

\begin{lemma}
\label{flagapartment}
For a pair of two flags $\{W_k\}$ and $\{W'_{k'}\}$, there exists an apartment containing both of them, iff they are compatible in the sense that for all elements $W_j \in \{W_k\}$  there are elements $W'_{j'},W'_{l'}\in \{W'_{k'}\}$ such that $W'_{j'} \subset W_{j} \subset W'_{l'}$ and vice versa.
Compatibility defines an equivalence relation on the space of flags.
\end{lemma}

\begin{cor}
 If the flags $\{W_k\}$ and $\{W'_{k'}\}$ are in $Gr_x^{n,+}$ then all flags in the apartment containing them are in $Gr_x^{n,+}$.
\end{cor}

\begin{proof}
 This is a consequence of lemma \ref{regularity conditions for frames}, condition $(iii)$.
\end{proof}

\begin{cor}
A flag $\{W_k\}$ is contained in an apartment containing the standard normal flag iff $W_k \subset Gr_0^{n,+}$.
\end{cor}

\begin{proof}[Proof of lemma \ref{flagapartment}]
By lemma~\ref{chambercomplex}, the complex of all flags is a chamber complex. Hence without loss of generality we can assume that $\{W_k\}$ and $\{W_k'\}$ are two maximal compatible flags. For each $k \in \mathbb{Z}$, we define the set $\pi(k):= 
\{m|\exists v \in (W_m\backslash W_{m+1})\cap (W_k'\backslash W_{k+1}')\}$. We have to show that $|\pi(k)|=1$ for all $k$.
So for $i\in \{0,\dots, n-1\}$ we choose vectors $v_i \subset \pi(i)$ and put $U_i = \textrm{span}\langle v_i\rangle$. Furthermore, for $i'=ln+i$ set $U_{i'}=U_{ln+i}= z^l U_{i}$.

The proof now consists of several steps:
\begin{itemize}
\item[-]  $\{U_k\}$ is a periodic frame.  As the flags $\{W_k\}$ and $\{W_k'\}$ are periodic, the conditions $v\in (W_m\backslash W_{m+1})\cap (W_k'\backslash W_{k+1}')$ and $z^l v \in (W_{m+ln}\backslash W_{m+1+ln})\cap (W_{k+ln}'\backslash W_{k+1+ln}')$ are equivalent. 
\item[-]  $W_k= W_{k+1}\oplus U_k$, $W_m= W_{m+1}\oplus U_k$ for all $k$ and $m\subset \pi(k)$. So the apartment associated to $\{U_k\}$ contains $\{W_k\}$ and $\{W_k'\}$. 
\item[-]  $\pi(k+n)=\pi(k)+n$ follows from the periodicity of $\{W_k\}$ and $\{W_k'\}$.
\item[-]  So we are left with showing that $\pi$ is a permutation, that is $|\pi(k)|=1\hspace{1ex}\forall k$. 
The compatibility condition gives 
\begin{align*}z^{l+1}\{W_k\}&=\{W_{k+(l+1)n}\}\subset z\{W_k'\}\\
&=\{W_{k+n}'\}\subset\{W_{k+1}'\} \subset \{W_k'\}\subset z^{-l}\{W_k\}\,. \end{align*}

So $W'_k\backslash W'_{k+1} \subset W_{k-ln}$. This shows that there are numbers $m$ such that the set $(W_m\cap W_k'\backslash W_{k+1}')$ is nonempty. 
On the other hand $W_{k+(l+1)n} \subset W'_{k+1}$ shows that the set of those $m$ is bounded from above. 
So there is for every $k$ a maximal $m$ such that $(W_m\cap (W_k'\backslash W_{k+1}')$ is nonempty. 
But then $(W_{m+1}\cap (W_k'\backslash W_{k+1}'))$ is empty.
So $(W_m \backslash W_{m+1}')\cap (W_k'\backslash W_{k+1}')$ is nonempty. So $\pi(k)$ is nonempty for all $k$. Symmetrically also $\pi^{-1}(m)$ is nonempty for all $m$. We use now the periodicity condition: $\pi(n+k)=\pi(k)+n$.  As each set $\pi(k)$ is nonempty, for every $l\in \{0, \dots, n-1\}$ there is $k$, such that $\pi(k)=l(\textrm{mod } n)$. 

This means that $\pi(k+n\mathbb Z)=l+n\mathbb{Z}$ as for each $l$, $l+n \mathbb{Z}$ is in the image.  The pigeon hole principle asserts that $|\pi(k)|=1 (\textrm{mod n})$.
 Then the periodicity shows that $|\pi(k)|=1$. Hence $\pi$ is a permutation and thus an element of $W_{\textrm{aff}}$. 
\qedhere
\end{itemize}
\end{proof}

Next we investigate isomorphisms between apartments:

\begin{lemma} [apartments are isomorphic]
\label{axiombuilding2foralgebraicflagbuildings}

For every pair of apartments $A$ and $A'$, there is an isomorphism $\varphi : A\longrightarrow A'$. If the intersection $A\cap A'$ is nonempty, then one can choose $\varphi$ in a way that it fixes $\{W_k\}$ and $\{W'_{k'}\}$ pointwise.
\end{lemma}

\begin{proof}
Isomorphisms between apartments correspond to isomorphisms of their frames. So let $A$ and $A'$ denote two apartments with frames $\{U_k\}$ and $\{U'_l\}$. Then every bijective map  $\varphi:\{U_k\}\longrightarrow\{U'_l\}$ that preserves the periodicity condition, induces an isomorphism of the apartments.
 If there is a full flag $\{W_k\} \subset A\cap A'$, then without loss of generality we can assume that $W_k:=\bigoplus_{l\geq k} U_l$ and $W_k:=\bigoplus_{l'\geq k} U'_l$. Then every map $\varphi_m: U_k \longrightarrow U'_{k+m}$ for all $m\in \mathbb{N}$ induces an isomorphism $\varphi: A\longrightarrow A'$, stabilizing $\{W_k\}$. If $\{W_k\}$ is a partial flag, then the preservation of the $\mathcal{K}$\ndash type restricts the possibilities for $m$. $\varphi_m$ preserves the cells in the boundary of $\{W_k\}$ iff $m\in n\mathbb{Z}$.

If $A\cap A'$ is nonempty and $\{W'_k\}\subset A\cap A'$, then this means just that \begin{displaymath}
W'_k:=\bigoplus_{l'\geq k} U'_{\pi(l)}=\bigoplus_{l'\geq k} U_{\pi(l)}. \, .
                                                                                    \end{displaymath}
Then $\phi_0$ trivially fixes the intersection. 
\end{proof}

\begin{theorem}[Tits building]
\label{titsbuilding}
Each equivalence class of compatible flags is an affine Tits building of type $A_n$.
\end{theorem}

\begin{cor}
The simplicial complex of positive (partial) flags in $Gr_0^{n,+}$ is an algebraic affine Tits building, called $\mathfrak{B}^+_0$.
\end{cor}

The proof of theorem~\ref{titsbuilding} is a consequence of lemmas~\ref{flagapartment} and~\ref{axiombuilding2foralgebraicflagbuildings} and corollary~\ref{chambercomplex_cor}:

\begin{proof}
We have to show that each equivalence class of compatible flags satisfies the axioms for a building: by corollary \ref{chambercomplex_cor} each equivalence class is a chamber complexes. The set of apartments consists of the set of all frames subjacent to full periodic flags.  The first axiom is lemma~\ref{flagapartment}, the second axiom is contained in lemma~\ref{axiombuilding2foralgebraicflagbuildings}.

Hence one finds an uncountable family of buildings which are all isomorphic to the well known building associated to the algebraic Kac-Moody group.
\end{proof}

As in the case of the flag manifolds the same constructions can be performed for the negative flag manifolds. In this way one obtains a second set of buildings, denoted $\mathcal{B^-}$. In section \ref{The generalized twin building} we show that $\mathfrak{B}=\mathfrak{B}^+ \cup \mathfrak{B}^-$ is a twin city.

This complex consists of all flags whose subspaces are in $Gr^{n,\pm}$. Those Grassmannians correspond to the loop group (resp.\ Kac-Moody groups) of $\frac{1}{2}$\ndash differentiable loops~\cite{PressleySegal86}. To get restrictions associated to smaller loop groups (i.e.\ loop and Kac-Moody groups of holomorphic loops) one needs the same constructions but applied to flags whose  subspaces are in the described sub-Grassmannians. Hence we get a family of twin cities corresponding to the different regularity conditions.

For the algebraic case this construction coincides with the well-known lattice description~\cite{Garrett97,AbramenkoNebe02,Kramer02}.

\subsection{The twin city}
\label{The generalized twin building}

To define a twinning between $\mathcal{B}^+$ and $\mathcal{B}^-$, we focus on the use of theorem~\ref{criterionforatwinbuilding}. Hence, we have to define a nonempty symmetric relation $\mathcal{O}$.
This is done via opposite flags:

\begin{definition}[Opposite flags]

\label{oppositenessrelation}
Let $\{W^+_k\}$ be a positive chamber and $\{W^-_k\}$ be a negative chamber. $\{W^+_k\}$ is opposite to $\{W^-_k\}$ iff $\textrm{dim } W^+_k \cap W^-_{n-k-1}=1$. 
\end{definition}

Using this definition, we define the symmetric relation $\mathcal{O}$ via the condition: $$(\{W^+_k\},\{V^-_k\})\in \mathcal{O}\textrm{ iff }\{W^+_k\}\textrm{  is opposite to }\{V^-_k\}\,.$$

This oppositeness relation induces a uniquely defined codistance function $d^*$ on $\mathcal{B}^+ \times \mathcal{B}^- \cup \mathcal{B}^- \times \mathcal{B}^+$.

\begin{theorem}
\label{hilberttwinbuildings}
The triple $(\mathcal{B}^+, \mathcal{B}^-, d^*)$ is a twin building.
\end{theorem}

To show that this twinning defines a twin building we verify the criterion of Ronan and van Maldeghem, (compare theorem~\ref{criterionforatwinbuilding}); this verification is contained in lemma~\ref{1twinning} and lemma \ref{criterionronanmaldeghem}:

\begin{lemma}
\label{1twinning}
$\mathcal{O}$ defines a $1$\ndash twinning.
\end{lemma}

\begin{proof}
Let $(\{W^+_k\},\{V^-_k\}), k\in \mathbb{Z}$ be a pair of opposite chambers, and let $(\{W'^+_k\},\{V'^-_k\})$ be walls of type $s_i$. Hence $(\{W'^+_k\},\{V'^-_k\}), k\in \mathcal{K}_{\hat{i}}$ are of type $m_{\mathcal{K}}= \{1,2, \dots, \hat{i}, \dots, n\}$. Hence at the $i$\ndash th position the frame $\mathcal{U}$ such that $(\{W'^+_k\}$ and $\{V'^-_k\})$ are subjacent to $\mathcal{U}$ is defined only up to a subspace $V^2$ of dimension $2$. The possible frames containing this wall correspond to the space of $1$\ndash dimensional subspaces in $V^2$ hence to $\mathbb{C}P_2$. Let $(e_1, e_2)$ denote a base of $V_2$ such that $\textrm{span}(e_1)=V_2 \cap W^+_i$. Then all chambers $\{V^-_k\}$ are opposite but the one defined by $V_2 \cap \{V^-_k\}=\textrm{span}(e_1)$.
\end{proof}

\begin{lemma}
\label{criterionronanmaldeghem}
For $\epsilon \in \{+,-\}$, there exists a chamber $c_{-\epsilon} \in \mathcal{C}_{-\epsilon}$ such that for any chamber $x_{\epsilon}$ with $(c_{-\epsilon}, x_{\epsilon})\subset \mathcal{O}$ there is an apartment $\Sigma_{\epsilon}$ of $\mathfrak{B}_{\epsilon}$ which satisfies $\{x_{\epsilon}\}= \{y_{\epsilon}\in \mathcal{C}_{\epsilon}|y_{\epsilon} \in \Sigma_{\epsilon} \textrm{ and } (c_{-\epsilon}, y_{\epsilon} ) \in \mathcal{O}\}$.
\end{lemma}

This condition follows directly from the description of apartments via frames: A frame gives rise to two apartments, one in $\mathfrak{B}^+$ and one in $\mathfrak{B}^-$.  This trivially satisfies the required condition. 

So we have proved the existence of a twinning between every pair of a positive and a negative building.

We call a city symmetric~\cite{Freyn10d}, iff there is an involution $\varphi$ such that $\varphi: \mathcal{B^{\pm}}=\mathcal{B}^{\mp}$
The symmetry between the positive and the negative Grassmannian yields the result:

\begin{cor}
The twin city $(\mathcal{B}^+, \mathcal{B}^-, d^*)$ is symmetric.
\end{cor}

\section[The symplectic groups \protect\boldmath$\widetilde{C}_{\textrm{\MakeLowercase{$n$}}}$]{The symplectic groups \protect\boldmath$\widetilde{C}_{n}$}
\label{The_symplectic_groups}

\subsection{The finite dimensional blueprint: type $C_n$}

The complex simple Lie group of type $C_n$ is the symplectic group $Sp(2n, \mathbb{C})$.
 $Sp(2n)$ acts on $V^{2n}:=\mathbb{C}^{2n}$. Let $\{e_1, \dots, e_n, f_1, \dots, f_n\}$ be a basis of $V^{2n}$ and define a symplectic form  $\omega$ by 
 $\omega( e_i, f_j) =\delta_{ij}= - \omega(f_j, e_i)$ and $\omega( e_i, e_j)=\omega( f_i, f_j )=0$. $Sp(2n)$ is the group of linear transformations preserving the symplectic form.

The Dynkin diagram of type $C_n$ is 

\centerline{
\begin{picture}(140,20)
\multiput(36,9)(30,0){5}{\circle{5}} 
 \put(33,0){\small $\alpha_1$\normalsize}\put(63,0){\small $\alpha_2$\normalsize}\put(91,0){\small $\alpha_{n-2}$\normalsize}\put(121,0){\small $\alpha_{n-1}$\normalsize}\put(153,0){\small $\alpha_{n}$\normalsize}
\multiput(42,9)(60,0){2}{\line(22,0){16}}\put(132,6){$\Longleftarrow$}
\put(74.5,8){$\dots$}
\end{picture} 
}

\noindent The Weyl group can be realized as a permutation group of any symplectic basis 
\begin{displaymath}
\{e_1, \dots, e_n, f_1, \dots, f_n\}\, .
\end{displaymath}
\noindent In this representation its generators $\alpha_k,k=1\dots n-1$ act as follows 
\begin{align*}
\alpha_k(e_k)&=e_{k+1}\qquad	&\alpha_k(f_k)&=f_{k+1}\\
\alpha_k(e_{k+1})&=e_{k+1}	&\alpha_k(f_k)&=f_{k}\\
\alpha_k(e_l)&=e_{l}		&\alpha_k(f_l)&=f_{l}, l\not=k,k+1\, .
\end{align*}

\noindent $\alpha_n$ acts as

\begin{align*}
\alpha_n(e_n)&=f_{n}\qquad	&\alpha_n(f_n)&=e_{n}\\
\alpha_n(e_l)&=e_{l}		&\alpha_k(f_l)&=f_{l}, l\not=n\,.
\end{align*}

To describe the spherical building of type $C_n$ we use the flag complex of isotropic subspaces:

\begin{definition}
A subspace $V\subset \mathbb{C}^n$ is called isotropic if $V\subset V^{\perp}$. A nested sequence of isotropic subspaces $V_1 \subset V_2 \subset \dots \subset V_l$ is called an isotropic flag. 
\end{definition}

An isotropic flag $F$ in $V^{2n}$ has length $l(V)\leq n$. Let $\{V_1, \dots, V_k\}$ be an isotropic flag. There is a unique completion 
\begin{displaymath}
\bar{F}:=\{V_1, \dots, V_k, V_{k-1}^{\perp}, \dots, V_1^{\perp}, V^{2n}\}\, .
\end{displaymath}
 In case of no ambiguity, we call this object ``isotropic flag''.

\begin{definition} 
Define $\Delta_{Sp(n)}$ to be the poset of isotropic flags with a partial order relation defined by inclusion. 
\end{definition}

\begin{theorem}
$\Delta_{Sp(n)}$ is a thick spherical building of type $C_n$.
\end{theorem}

\begin{proof}
\cite{AbramenkoBrown08}, section 6.6 or~\cite{Garrett97}, chapter 10.
\end{proof}

Let us mention a last construction: For an ascending isotropic flag $V_1\subset V_2 \subset \dots \subset V_n$ the sequence of their orthogonal complements $V_n^{\perp}\subset V_{n-1}^{\perp}\subset \dots \subset V_1^{\perp}$ form a second ascending flag. The cocatenation
\begin{displaymath}
V_1\subset V_2 \subset \dots \subset V_n = V_n^{\perp}\subset V_{n-1}^{\perp}\subset \dots \subset V_1^{\perp} 
\end{displaymath}

\noindent is a full flag. As each isotropic flag has a unique continuation in this way to a full flag, we can equivalently describe the building of type $C_n$ as the simplicial complex associated to the flag complex of completed flags. Remark that the symplectic group preserves orthogonality relations as it preserves the symplectic form; hence it maps completed flags onto completed flags.

\subsection{Symplectic Grassmannians}

The loop group $L(Sp(2n, \textrm{Id}))$ consists of periodic loops into the symplectic group $Sp(2n)$.
The Kac-Moody group $\widehat{L}(Sp(2n),\textrm{Id})$ is the torus extension over this loop group~\cite{PressleySegal86,Popescu05,Freyn09}.

In the finite dimensional blueprint the building is defined to be the flag complex over the set of all isotropic flags, ordered by inclusion. Equivalently there is a unique completion of flags using coisotropic subspaces. The same construction works in the infinite dimensional case. 

Let $H^{2n}:=L^2(S^1, V^{2n})$ and let $\left\{ z^{k}e_i,z^{l}f_{j}\right\}$ for $\{i,j\}\in \{1, \dots, n\}$ and  $\{k,l\}\in \mathbb{Z}$ be a basis of $H^{2n}$. The symplectic form $\Omega$ on $H^{2n}$ is defined by the condition 

\begin{displaymath}
\Omega( z^{k}e_{i},z^{l} f_{j}) =\delta_{i,j}\delta_{k,-l}= - \Omega (z^{l} f_{j}, z^{k} e_{i})\,.
 \end{displaymath}

The restriction of $\Omega$ to the $2n$\ndash dimensional subspace of constant loops is the nondegenerate $2$\ndash form $\omega$.
In geometric algebra a ``hyperbolic plane'' in a symplectic space is a $2$\ndash dimensional subspace $V=\textrm{span}\langle v_1, v_2\rangle $ such that $\omega(v_1, v_2)=-\omega(v_2, v_1)=1$ and  $\omega(v_1, v_1)=\omega (v_2, v_2)=0$. $(v_1, v_2)$ is called a hyperbolic pair~\cite{Garrett97}. 
A hyperbolic space is subjacent to a frame $\{U_k\}$ if it admits a hyperbolic pair that consists of subspaces $U_{k_1},U_{k_2}\in \{U_k\}$ A frame in a symplectic vector space is called symplectic if its elements can be paired into hyperbolic pairs. 
Hyperbolic planes in $H^{2n}$ are for example $V_{i,l}=\textrm{span}\langle z^le_i, z^{-l}f_i\rangle$ for $l\in \mathbb{Z}$ and $i\in \{1, \dots, n\}$.
Conversely all hyperbolic planes subjacent to the frame described by $\left\{ z^{k}e_i,z^{l}f_{j}\right\}$ for $\{i,j\}\in \{1, \dots, n\}$ and  $\{k,l\}\in \mathbb{Z}$.

\begin{definition}[Symplectic Grassmannian]
\label{Symplectic_Grassmannian}
Let $X$ denote any regularity condition we introduced in section \ref{Grassmanians_and_periodic_flags} and let $\epsilon\in \{+,-\}$.
\begin{enumerate}
\item The isotropic (resp.\ coisotropic, symplectic) Grassmannian of regularity $X$ denoted $Gr_{\textrm{iso},X}^{\epsilon}(H^{(2n)})$ (resp.\ $Gr_{\textrm{coiso},X}^{\epsilon}(H^{(2n)})$,  $Gr_{\textrm{symp},X}^{\epsilon}(H^{(2n)})$) consists of all isotropic (resp.\ coisotropic, isotropic and coisotropic) subspaces  in $Gr_X^{\epsilon}(H^{(2n)})$.
\item The reduced isotropic (resp.\ coisotropic, symplectic) Grassmannian of regularity $X$ denoted $Gr_{\textrm{iso},X}^{(2n), \epsilon}(H^{2n})$ (resp.\ $Gr_{\textrm{coiso},X}^{(2n), \epsilon}(H^{2n})$,  $Gr_{\textrm{symp},X}^{(2n), \epsilon}(H^{2n})$) consist of all isotropic (resp.\ coisotropic, isotropic and coisotropic) subspaces in $Gr_X^{(2n), \epsilon}(H^{2n})$.
\end{enumerate}
\end{definition}

For example $H$ is not isotropic hence not in the isotropic Grassmannian, but $zH\oplus K$ is if $K$ is an isotropic subspace in $V^{2n}$. Generally all  isotropic subspaces have virtual dimension less or equal than $-n$ while coisotropic subspaces have virtual dimension greater or equal than $-n$

As in the finite dimensional case each subspace has its orthogonal complement $V^{\perp}$. For example we have $H^{+}=(zH)^{\perp}\subset H^{+}$. The orthogonal complement of isotropic subspaces is coisotropic and vice versa.

\begin{remark}
Remark a crucial difference: In case $\widetilde{A}_n$ all subspaces are on an equal footing. There is no distinguished virtual dimension. Hence flags are stable under shifts. This is in crucial contrast to the situation of the isometry groups. Here exists in every flag a maximal {\it isotropic} subspace. Hence there are subspaces of distinguished virtual dimension. This functional analytic fact reflects the antipode between the cyclic Dynkin diagram of type $\widetilde{A}_n$ and the linear Dynkin diagram for the isometry groups.
\end{remark}

\begin{definition}[positive symplectic periodic flag manifold]

The full positive periodic symplectic flag manifold $Fl_{\textrm{symp}}^{+,n}$ consists of all sequences $W_k, k \in \mathbb{Z}$ of subspaces in $H^{n}$ such that 
\begin{enumerate}
	\item $W_k \in Gr_{\textrm{symp},X}^{\pm}(H^{(2n)})$ ,
	\item $W_{k+1} \subset W_k \forall k$ and $\textrm{dim}(W_k/W_{k+1})=1$, 
	\item $zW_k = W_{k+2n}$.
	\item $W_k=W_{2n-k}^{\perp}$
\end{enumerate}
\end{definition}

The definition of partial periodic submanifolds and their negative counterparts follows the pattern described for the case $A_n$.

Let us remark that~\cite{PressleySegal86} uses a different definition for ``symplectic'' Grassmannians:

\begin{displaymath}
Gr_{\mathbb{H}}^{(2n)}:=\{W\in Gr^{(2n)}|(jW)^{\perp}=zW\}\, .
\end{displaymath}

Here we omit the precise regularity assumptions. This definition is motivated by the equivalence $\Omega(Sp(n))\cong Gr_{\mathbb{H}}^{(2n)}$~\cite{PressleySegal86}, 8.5.
This Grassmannian consists of the component of virtual dimension $0$\ndash component of our (coisotropic) symplectic Grassmannian. As the action of the Kac-Moody group of type $\widetilde{C}^n$ preserves virtual dimension of subspaces the restriction to one virtual dimension is the appropriate choice for the investigation of the group action; nevertheless the choice of the component of codimension $-n$, consisting of Lagrangian subspaces, would perhaps have been more natural from a symplectic viewpoint but of course less natural from a functional analytic one.

We know~\cite{Freyn10d} from the abstract theory that a city of type $\widetilde{C}_n$ can be described by 
\begin{displaymath}
\mathfrak{B}:=\left(L(Sp(n), \textrm{Id})/B \times \Delta\right)/\sim\, .
\end{displaymath}
Using the equivalence $\Omega(Sp(n))\cong Gr_{\mathbb{H}}^{(2n)}$ it is obvious that this description can be made explicit by the careful construction of the appropriate flags of suitable isotropic (resp.\ symplectic) subspaces in $Gr_{\mathbb{H}}^{(2n)}$~\cite{Freyn09}.

\subsection{The affine building of type $\widetilde{C}_n$}

Let us state the main result:

\begin{theorem}
\label{twin_city_of_type_C_n}
The simplicial complex associated to the poset of symplectic periodic flags is a twin city of type $\widetilde{C}_n$. 
\end{theorem}

For the proof of this theorem we need a description of apartments; to this end we use the following class of frames:

\begin{definition}
 
A frame consists of elements $\left\{ e_{ik},f_{jl}, \{i,j\}\in \{1, \dots, n\}, \{k,l\}\in \mathbb{Z}\right\}$ such that
\begin{itemize}
  \item[-] $\left\{ e_{i,k},f_{j,l}, \{i,j\}\in \{1, \dots, n\}, \{k,l\}\in \mathbb{Z}\right\}$ span $H^{2n}$, 
  \item[-] $z e_{ik}=e_{i,k+1}$ and $\frac{1}{z} f_{jl}= f_{j,l-1}$,
  \item[-] $e_{i,k}$ and $f_{j,l}$ are isotropic lines,
  \item[-] $(e_{i,k}\oplus f_{i,-k})$ is a hyperbolic plane and $(e_{i,k}\oplus f_{i,-k})\perp (e_{j,l}\oplus f_{j,-l})$ for all $i\not=j\in \{1, \dots, n\}, k,l\in \mathbb{N}$.
\end{itemize}
\end{definition}

Remark that those frames are hyperbolic.

\begin{definition}
The apartment described by the frame $\left\{ e_{ik},f_{j,-l}, \{i,j\}\in \{1, \dots, n\}, \{k,l\}\in \mathbb{N}\right\}$ consists of all flags subjacent to this frame.
\end{definition}

The various regularity conditions for frames are defined in complete analogy to the case of $\widetilde{A}^n$.

\begin{definition}
A permutation of the frame $\left\{ e_{ik},f_{jl}, \{i,j\}\in \{1, \dots, n\}, \{k,l\}\in \mathbb{Z}\right\}$ is admissible if it preserves periodicity and hyperbolic planes.
\end{definition}

\begin{lemma}
The group of admissible permutations is a Weyl group of type $\widetilde{C_n}$. 
\end{lemma}

\begin{proof}
The Dynkin diagram of type $\widetilde{C}^n$ looks like that:

\centerline{
\begin{picture}(140,30)
\multiput(6,15)(30,0){6}{\circle{5}} 
\put(3,6){\small$\alpha_0$\normalsize} \put(33,6){\small $\alpha_1$\normalsize}\put(63,6){\small $\alpha_2$\normalsize}\put(91,6){\small $\alpha_{n-2}$\normalsize}\put(121,6){\small $\alpha_{n-1}$\normalsize}\put(153,6){\small $\alpha_{n}$\normalsize}
\multiput(42,15)(60,0){2}{\line(22,0){18}}\put(12,12){$\Longrightarrow$}\put(132,12){$\Longleftarrow$}
\put(74.5,14){$\dots$}
\end{picture} 
}

Take first the $2n$\ndash dimensional subspace $H= \textrm{span}\{e_i, f_i| i=1, \dots, n\}$. The Weyl group is of type $C_n$, corresponding to the Dynkin diagram defined by deleting $\alpha_0$. It is the semidirect product of the symmetric group in $n$\ndash letters transposing the indices and a $\mathbb{Z}_2$ group for each index, interchanging $e_i$ and $f_i$ and fixing all other bases. We are left with checking the action of the additional vertex $\alpha_0$. It acts by interchanging $e_1$ and $zf_1$. The check of the necessary commutation relations is straight forward. This vertex generates an additional $\mathbb{Z}^n$\ndash factor in $\widetilde{C}_n$, which acts by raising or lowering the powers of $z$. 
\end{proof}

To prove theorem \ref{twin_city_of_type_C_n} we have to verify several details. The steps are analogous to the case of $\widetilde{A}^n$.

\begin{enumerate}
	\item The apartments are thin Coxeter complexes of type $\widetilde{C}_n$.
	\item The whole complex is thick.
	\item Each connected component is a building, hence
	\begin{enumerate}
		\item  Every pair of chambers is contained in a common apartment. 
		\item All apartments are isomorphic.
	\end{enumerate}
	\item Check that the whole complex is a twin city.  
\end{enumerate}

\begin{lemma}
 The apartments are thin Coxeter complexes of type $\widetilde{C}_n$.
\end{lemma}

\begin{proof}
 Let $\mathcal{F}:=\{e_{ik},f_{jl}\}$ be a symplectic frame. An isotropic flag $\{W_k\}$ is subjacent to $\mathcal{F}$ if each subspace $W_k$ is spanned by a subset of elements in $\mathcal{F}$. Let $i\in \{1, \dots n\}$ and let $\{W_k\}$ be of type $\mathcal{K}=\{1,2, \dots i-1, i+1, \dots, 2n-i-1, 2n-i+1, \dots, 2n\}$. Flags of this type correspond to codimension $1$\ndash simplices. The subspace $W_{i+1}/W_{i-1}$ has dimension $2$; hence  modulo $W_{i-1}$ it is generated by two elements in $\mathcal{F}$. Thus there are two possible ways to construct a subspace $W_i$ subjacent to $\mathcal{F}$ such that $W_{i-1}\subset W_i\subset W_{i+1}$;  the same is true for $W_{2n-i+1}/W_{2n-i-1}$. The condition to be symplectic flags assures that a choice of $W_i$ defines unambiguously $W_{2n-i}$.  Hence there are exactly two possible completions of $\{W_k\}$ to full flags. Thus the apartment is thin. From the action of the Weyl group of type $\widetilde{C_n}$ on a frame we see that it acts transitively on the set of maximal flags - hence it is of type $\widetilde{C_n}$.
\end{proof}

\begin{lemma}
 The whole complex is thick.
\end{lemma}

\begin{proof}
We need to check that any codimension $1$\ndash simplex is in the boundary of at least $3$ chambers. For codimension $1$ \ndash simplices corresponding to  $\alpha_1, \dots \alpha_n$, this follows exactly as in the finite dimensional case, for $\alpha_0$ the situation is identical to the one of $\alpha_n$.
\end{proof}

Connected components are defined exactly in the same way as for $\widetilde{A}^n$.

\begin{lemma}
In a connected component, each pair of chambers is contained in a common apartment. 
\end{lemma}
\begin{proof}
Two flags describe cells in the same connected component, iff they are compatible in the sense that for all elements $W_j \in \{W_k\}$  there are elements $W'_{j'},W'_{l'}\in \{W'_{k'}\}$ such that $W'_{j'} \subset W_{j} \subset W'_{l'}$ and vice versa. Compatibility clearly defines an equivalence relation on the space of flags. For compatible flags the strategy of proof, described in~\cite{Garrett97}, chapter 10  and 20 can be adapted. 
\end{proof}

\begin{lemma}
All apartments are isomorphic.
\end{lemma}
\begin{proof}
Isomorphism of apartments correspond to isomorphisms of their frames - hence this axiom is clear. 
\end{proof}

\begin{lemma}
The whole complex is a twin city.
\end{lemma}

\begin{proof}
We define the codistance via opposite flags exactly as we did in case of $\widetilde{A}_n$, definition \ref{oppositenessrelation}. The verification is straight forward.
\end{proof}

\section[The orthogonal groups \protect\boldmath$\widetilde{B}_{\textrm{\MakeLowercase{$n$}}}$ and \protect\boldmath$\widetilde{B}_{\textrm{\MakeLowercase{$n$}}}$]{The orthogonal groups \protect\boldmath$\widetilde{B}_{n}$ and \protect\boldmath$\widetilde{D}_{n}$}
\label{The_orthogonal_groups}

\subsection{The finite dimensional blueprint: types \protect\boldmath$B_n$ and \protect\boldmath$D_n$}

Similar to the symplectic groups we start with an $m$\ndash dimensional vector space $V^{m}(\mathbb{K})$, $\mathbb{K} \in \{\mathbb{R}, \mathbb{C}\}$ and define a non-degenerate quadratic form $q$ on $V^{m}$ in the following way: For a standard basis $e_1, \dots, e_n$ of $V^{m}$ put 
\begin{displaymath}
q(e_i, e_j)=\delta_{i+j,m+1}\,.
\end{displaymath}

\begin{itemize}
\item[-] If $\mathbb{K}=\mathbb{C}$ then this quadratic form is equivalent to the standard quadratic form. Hence maximal isotropic subspaces have dimension $n$. If $\mathbb{K}=\mathbb{R}$ then the form has Witt index $n=\lfloor \frac{m}{2}\rfloor$ on $V^{m}$.  The orthogonal group $O_q(m, \mathbb{K})$ is the group of linear automorphisms of $V^{m}$ that preserve $q$.
We study the isotropic subspaces. 
If $m=2n+1$ then the group $O_q(m, \mathbb{K})$ is of type $B_n$; its Dynkin diagram is:

\centerline{
\begin{picture}(165,30)
\multiput(36,15)(30,0){5}{\circle{5}} 
 \put(33,6){\small $\alpha_1$\normalsize}\put(63,6){\small $\alpha_2$\normalsize}\put(91,6){\small $\alpha_{n-2}$\normalsize}\put(121,6){\small $\alpha_{n-1}$\normalsize}\put(153,6){\small $\alpha_{n}$\normalsize}
\multiput(42,15)(60,0){2}{\line(22,0){18}}\put(132,12){$\Longrightarrow$}
\put(74.5,14){$\dots$} 
\end{picture} 
}

The associated building can be constructed equivalently to the case of $C_n$ as the flag complex of subspaces isotropic with respect to $q$. This construction yields a building of type $B_n$ if the Witt index of $q$ is $n$. For details see~\cite{Garrett97} chapter 10.

\item[-] If $m=2n$ then $O_q(m, \mathbb{K})$ is a group of type $D_n$;
its Dynkin diagram is

\centerline{
\begin{picture}(160,60)
\multiput(36,21)(30,0){4}{\circle{5}} \put(147,0){\circle{5}}\put(147,42){\circle{5}}
 \put(33,27){\small $\alpha_1$\normalsize}\put(63,27){\small $\alpha_2$\normalsize}\put(91,27){\small $\alpha_{n-3}$\normalsize}\put(133,21){\small $\alpha_{n-2}$\normalsize}\put(148,5){\small $\alpha_{n-1}$\normalsize}\put(148,46){\small $\alpha_{n}$\normalsize}
\multiput(42,22)(60,0){2}{\line(22,0){18}}     \put(129,25){\line(1,1){15}}\put(129,17){\line(1,-1){15}}
\put(74.5,21){$\dots$}
\end{picture} 
}

 We have to distinguish two cases: if the group is not $O(n,n)$ then the flag complex of isotropic subspaces is a thick spherical building of type $C_n$. Contrary if the Witt index is $n$ then the complex of isotropic subspaces is not thick.  It turns out that the associated spherical building of type $D_n$ is equivalent to the oriflamme complex associated to $q$~\cite{AbramenkoBrown08}, section 6.7.\ and~\cite{Tits74,Garrett97}.

To construct the oriflamme complex we start with the set of all nonzero subspaces of $V^{2n}$, isotropic with respect to $q$ and define the building as the flag complex of a complex whose full flags are all of the form

\begin{displaymath}
V_1\subset V_2 \subset \dots \subset V_{n-2}\subset V_{n,1}, V_{n,2}
\end{displaymath}
such that $\textrm{dim}V_i=i$. Hence our flags contain no subspace of dimension $n-1$, but they contain two subspaces of dimension $n$. We put $V_{n-1}=V_{n,1}\cap V_{n,2}$ and require furthermore that $\textrm{dim}V_{n-1}=\textrm{dim} (V_{n,1}\cap V_{n,2})=n-1$. Clearly $V_{n-2}\subset V_{n-1}$. Frames consist of lines $\langle p_i\rangle$ and $\langle q_i\rangle$, $i\in \{1, \dots n\}$ such that for each $i$ the space $\textrm{span}\langle p_i, q_i \rangle$ is a hyperbolic plane.
\end{itemize}

\subsection{Orthogonal Grassmannians}
\label{orthogonal_grassmannians}

The constructions for $\widetilde{B}_n$ and $\widetilde{D}_n$ are similar. The problem is just that the complex of isotropic flags, as we defined it for the cities of type $\widetilde{C}_n$, is not thick. Hence we have to review our construction to get more refined combinations of nested subspaces to represent the simplices of the building. The simplicial complexes constructed in this way are called oriflamme complexes. For $\widetilde{D}_n$ we need the double oriflamme complex, for $\widetilde{B}_n$ the simple oriflamme complex. The analytic part of the theory proceeds similarly to the case of $\widetilde{A}_n$ and $\widetilde{C}_n$. For the algebraic details of the oriflamme construction compare~\cite{Garrett97}, chapter 20. Affine buildings of this type are constructed in~\cite{Garrett97,AbramenkoNebe02}.

The invariant form on $H^{2n}$ is defined by the condition 
\begin{displaymath}
Q(z^{k}e_{i},z^{l} f_{j}\rangle =q( e_{i}, f_{j}) \delta_{k+l,0}=  q( z^{l} f_{j}, z^{k} e_{i} )\,.
\end{displaymath}

\noindent where $\{e_i, f_i\}$ is a basis for $H^0$ and $q$ is a scalar product on $V^{2n}$ of Witt index $n$.

\begin{definition}[Orthogonal Grassmannian]
\label{orthogonal_grassmannian}
Let $X$ denote any of the regularity conditions we introduced in section \ref{Grassmanians_and_periodic_flags} and $\epsilon\in \{+,-\}$.
\begin{enumerate}
\item The isotropic (resp.\ coisotropic, orthogonal) Grassmannian of regularity $X$ denoted $Gr_{\textrm{iso},X}^{\epsilon}(H^{(m)})$ (resp.\ $Gr_{\textrm{coiso},X}^{\epsilon}(H^{(m)})$,  $Gr_{\textrm{orth},X}^{\epsilon}(H^{(m)})$) consists of all isotropic (resp.\ coisotropic, isotropic and coisotropic) subspaces  in $Gr_X^{\epsilon}(H^{(m)})$.
\item The reduced isotropic (resp.\ coisotropic, symplectic) Grassmannian of regularity $X$ denoted $Gr_{\textrm{iso},X}^{(m), \epsilon}(H^{m})$ (resp.\ $Gr_{\textrm{coiso},X}^{(m), \epsilon}(H^{m})$,  $Gr_{\textrm{orth},X}^{(m), \epsilon}(H^{m})$) consist of all isotropic (resp.\ coisotropic, isotropic and coisotropic) subspaces in $Gr_X^{(m), \epsilon}(H^{m})$.
\end{enumerate}
\end{definition}

\noindent Let us remark that~\cite{PressleySegal86}, section 8.5. (not going into the details of regularity questions) uses the definition 
\begin{displaymath} Gr^{(n)}=\{W\in Gr^{(n)}|\overline{W}^{\perp}=zW\}\, .\end{displaymath}
This corresponds to the connected component of virtual dimension $0$ of our Grassmannian. The justification follows similar to the case of symplectic groups from the correspondence between $Gr^{(n)}$ and $\Omega (O^n)$.

\subsection{The double oriflamme complex}

We start with a metric as defined in section \ref{orthogonal_grassmannians}. If we are working over $\mathbb{R}$ suppose that the restriction of the metric to $V^{2n}$ has Witt index $n$.

\begin{definition}[positive periodic double oriflamme manifold]

The full positive periodic orthogonal oriflamme manifold $Fl_{\textrm{ort}}^{+,n}$ consists of all sets $W_k, k \in \mathbb{Z}$ of subspaces
 in $H^{n}$ such that 

\begin{enumerate}
	\item $W_k \subset Gr_{\textrm{ort},X}^{\pm, \textrm{ext}}(H^{(2n)})$ and $\textrm{virt dim}(W_k)(\textrm{mod }n)\not= \pm 1$,
	\item If $k(\textrm{mod }n)=0$, then there are two subspaces $W_{k,1}$ and $W_{k,2}$ such that 
$$\textrm{virt dim }\left(W_{k,1}\cap W_{k,2}\right)(\textrm{mod }n)=-1\, .$$ 
Denote this space $W_{k+1}$. Furthermore we put $W_{k-1}=W_{k,1}\cup W_{k,2}$,.
	\item $W_{k+1} \subset W_k$ and $\textrm{dim}(W_k/W_{k+1})=1$. 
	\item $zW_k = W_{k+n}$.
\end{enumerate}
\end{definition}

\noindent Hence a double oriflamme consists of a series of subspaces 
\begin{align*}
\dots \subset W_{0,1}, W_{0,2}\subset W_{0,1}\cup W_{0,2}&\subset W_{-2}\subset W_{-3} \subset \dots\\
& \dots \subset W_{-n+2}\subset W_{-n,1}, W_{-n,2}\subset W_{-n-2}\subset \dots
\end{align*}

\noindent The definition of positive partial periodic double oriflammes is self explaining.

\begin{theorem}[Double oriflamme complex]
\label{double_oriflame}
 The simplicial complex associated to the poset of positive and negative partial periodic single oriflammes is an affine twin city of type $\widetilde{D}_n$. 
\end{theorem}

Apartments are defined via frames:

\begin{definition}
A frame consists of a set of lines $\langle p_{i,k}\rangle, \langle q_{i,l}\rangle, i=1,\dots, n$ and ,$k,l\in \mathbb{Z} $ such that 
\begin{enumerate}
 \item $H_{i,k} =\textrm{span}\langle p_{i}, q_{i}\rangle$ is a hyperbolic plane and $\langle p_i\rangle, \langle q_i\rangle$ is a hyperbolic pair.
 \item For $i\not=j$ we have $H_i \perp H_j$.
 \item $z \langle p_{i,k} \rangle=\langle p_{i,k+1}\rangle$ and $z\langle q_{i,k}\rangle=\langle q_{i,k+1}\rangle$.
 \item The $z$\ndash cyclic subspace of spanned by $\langle p_{i,k_i}\rangle, \langle q_{i,l_i}\rangle, i=1, \dots, n$ is in the orthogonal Grassmannian $Gr_{\textrm{ort},X}^{\pm, \textrm{ext}}(H^{(2n)})\, $.
\end{enumerate}
\end{definition}

\begin{proof}
The proof of this theorem follows the case $\widetilde{C}_n$. We need the following facts (compare~\cite{Garrett97}):

\begin{enumerate}
	\item The apartments are thin Coxeter complexes of type $\widetilde{D}_n$.
	\item The whole complex is thick.
	\item Each connected component is a building. We need to know:
	\begin{enumerate}
		\item  Every pair of chambers is contained in a common apartment. 
		\item All apartments are isomorphic.
	\end{enumerate}
	\item Check that the whole complex is a twin city (hence check the twinning).
\end{enumerate}
\end{proof}

To prove that apartments are thin we have to check that each partial periodic flag omitting exactly one subspace in each period has exactly two extensions. This is a straight forward verification of possible cases. The proof of thickness of the whole complex is similar. The double oriflamme complex is associated to groups of type $\widetilde{D}^n$. The Dynkin diagram is

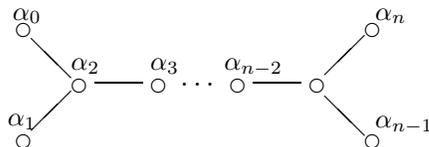
\centerline{
\begin{picture}(160,60)
\multiput(36,21)(30,0){4}{\circle{5}} \put(15,0){\circle{5}}\put(15,42){\circle{5}}\put(147,0){\circle{5}}\put(147,42){\circle{5}}
 \put(11,46){\small $\alpha_{0}$\normalsize}\put(9,6){\small $\alpha_{1}$\normalsize} \put(33,27){\small $\alpha_2$\normalsize}\put(63,27){\small $\alpha_3$\normalsize}\put(91,27){\small $\alpha_{n-2}$\normalsize}\put(148,5){\small $\alpha_{n-1}$\normalsize}\put(148,46){\small $\alpha_{n}$\normalsize}
\multiput(42,22)(60,0){2}{\line(22,0){18}} \put(18,4){\line(1,1){15}}\put(18,39){\line(1,-1){15}}\put(129,25){\line(1,1){15}}\put(129,17){\line(1,-1){15}}
\put(74.5,21){$\dots$} 
\end{picture} 
}

To prove that apartments are of type $\widetilde{D}_n$ we have to check that the Weyl group acts transitively on the set of all full oriflammes subjacent to a frame. This is done by checking the action of every generator of the Weyl group $\alpha_i$. The first axiom: Every pair of chambers is contained in a common apartment follows from an adaption of the argument used in the case $\widetilde{A}_n$;
the isomorphisms between two arbitrary apartments are constructed via isomorphisms of frames. To prove that the whole complex is a twin building we can check any of the two conditions given in appendix \ref{some_results_about_twin_buildings}.

\subsection{The single oriflamme complex}

The single periodic oriflamme complex corresponds to the buildings of type $\widetilde{B}_n$. Groups of type $\widetilde{B}_n$ are the Kac-Moody groups associated to orthogonal groups of odd dimension. Hence let $H^n:=L^2(S^1, \mathbb{F}^{2n+1})$ with $\mathbb{F}\in \{\mathbb{R}, \mathbb{C}\}$ and define a scalar product $\langle \cdot, \cdot\rangle$ by

\begin{displaymath}
\langle(q_{i,1}, p_{i,1}, r_1),(q_{j,2}, p_{j,2}, r_2)\rangle:=\sum q_{i,1}p_{i,2}+\sum q_{i,2}p_{i,1}+ \eta r_1 r_2\, .
\end{displaymath}

\begin{definition}[positive periodic single oriflamme manifold]

The full positive periodic single oriflamme manifold $Fl_{\textrm{ort}}^{+,n}$ consists of all sets of subspaces $W_{k,j}, k \in \mathbb{Z}, j\in \{1,2\}$  in $H^{n}$ such that $k (\textrm{mod } 2n)\not= \pm 1,  j=1$ if $k(\textrm{mod } 2n)\not= 0$ 
satisfying the following conditions:

\begin{enumerate}
	\item $W_{k,j} \subset Gr_{\textrm{ort},X}^{\pm, \textrm{ext}}(H^{(2n)})$ and $\textrm{virt dim}(W_{k,j})= -k$.
	\item If $k(\textrm{mod }2n)=0$, then there are two subspaces $W_{k,1}$ and $W_{k,2}$ such that 
$$\textrm{virt dim }\left(W_{k,1}\cap W_{k,2}\right)(\textrm{mod }2n)=-1\, .$$ 
	\item $W_{k,j} \subset W_{l,j}$ and $\textrm{dim}(W_k/W_{k+1})=1$ if $k > l$.
	\item $zW_k = W_{k+2n}$.
\end{enumerate}
\end{definition}

Hence a single oriflamme consists of a series of subspaces 
\small\begin{align*} 
\dots &\subset W_{0,1}\cap W_{0,2}\subset W_{0,1}, W_{0,2}\subset W_{0,1}\cup W_{0,2}\subset W_{-2}\subset W_{-3} \subset \dots\\
&\dots \subset W_{-n+2}\subset W_{-n+1} \subset W_{-n}\subset W_{-n-1}\subset W_{-n-2}\subset \dots 
\end{align*}\normalsize

\noindent The definition of positive partial periodic single oriflammes is self explaining.

\begin{theorem}
 The simplicial complex associated to the poset of positive and negative partial periodic single oriflammes is an affine twin city of type $\widetilde{B}_n$. 
\end{theorem}

\noindent Apartment are defined via frames:

\begin{definition}
A frame consists of a set of lines $\langle p_{i,k}\rangle, \langle q_{i,k}\rangle, i=1, \dots n, k\in \mathbb{Z}$ such that 
\begin{enumerate}
 \item $H_i =\textrm{span}\{ \langle p_i\rangle, \langle q_i\rangle\}$ is a hyperbolic plane and $\langle p_i\rangle, \langle q_i\rangle$ is a hyperbolic pair.
 \item For $i\not=j$ we have $H_i \perp H_j$.
 \item The $z$\ndash cyclic subspace of spanned by $\langle p_{i,l_i}\rangle, \langle q_{i,k_i}\rangle, i=1, \dots n$ is in the orthogonal Grassmannian $Gr_{\textrm{ort},X}^{\pm, \textrm{ext}}(H^{(2n)})$.
\end{enumerate}

\end{definition}

The rest is straight forward verification along the lines of the proof of theorem \ref{double_oriflame}. The main difference is that we have only one oriflamme and that apartments are now of type $\widetilde{B}_n$. The Dynkin diagram is 

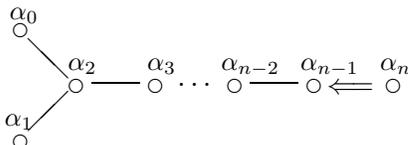
\centerline{
\begin{picture}(160,60)
\multiput(36,21)(30,0){5}{\circle{5}} \put(15,0){\circle{5}}\put(15,42){\circle{5}}
 \put(11,46){\small $\alpha_{0}$\normalsize}\put(9,6){\small $\alpha_{1}$\normalsize} \put(33,27){\small $\alpha_2$\normalsize}\put(63,27){\small $\alpha_3$\normalsize}\put(91,27){\small $\alpha_{n-2}$\normalsize}\put(121,27){\small $\alpha_{n-1}$\normalsize}\put(151,27){\small $\alpha_{n}$\normalsize}
\multiput(42,22)(60,0){2}{\line(22,0){18}} \put(18,4){\line(1,1){15}}\put(18,39){\line(1,-1){15}}
\put(74.5,21){$\dots$} \put(131,18){$\Longleftarrow$}
\end{picture} 
}

Thus it is a combination of one end from the case $\widetilde{C}_n$ and one end from the case $\widetilde{D}_n$.

\appendix
\section{Some results about twin buildings}
\label{some_results_about_twin_buildings}

In this appendix we recall some characterizations of buildings and twin buildings, we use in our proofs. Main references are~\cite{Garrett97,AbramenkoBrown08}.
 
We start with the definition of buildings and twin buildings:

\begin{definition}
 A building is a thick chamber complex $\Delta$, together with a set of thin subcomplexes, called apartments $\mathcal{A}$ such that
\begin{enumerate}
 \item For each pair of chambers $x,y\in \Delta$ there is $A_{xy}\in \mathcal{A}$ containing both.
 \item For any pair of two apartments $A, A'\in \mathcal{A}$, there is a chamber complex isomorphism: $\varphi: A\longrightarrow A'$ fixing the intersection $A\cap A'$ 
\end{enumerate}
\end{definition}

A chamber complex is a connected simplicial complex such that each cell is contained in a cell of maximal dimensions and for each pair of chambers $x$ and $y$ there is a sequence of chambers $x_i$, such that $x_1=x$ and $x_n=y$ and $x_i\cap x_j$ is a wall. A chamber complex is thick if each wall is in the boundary of at least three chambers, it is thin if each wall is in the boundary of $2$ chambers.

\begin{definition}[Twin building]
\label{twinbuilding}
A twin building consists of a pair of buildings $\mathfrak{B}^+$ and $\mathfrak{B}^-$ together with a codistance function $\delta^*_W: \mathfrak{B}^+\times \mathfrak{B}^-\longrightarrow W$ which is  subject to the following conditions: Let $X \in \mathfrak{B}^+$ and $Y, Z \in \mathfrak{B}^-$,
\begin{enumerate}
\item $\delta^*(X,Y)= \delta^*(Y,X)^{-1}$,
\item $\delta^*(X,Y)=w, \delta(Y,Z)=s \in S$ and $l(ws)=l(w)-1$, then $\delta^*(X,Z)=ws$,
\item $\delta^*(X,Y)=w, s \in S$. Then $\exists Z\in \mathfrak{B}^-$ such that $\delta^*(Y,Z)=s$ and $\delta^*(X,Z)=ws$.
\end{enumerate}
\end{definition}

Two common characterization of twin buildings are the following: The first one due to Peter Abramenko and Hendric Van Maldeghem~\cite{AbramenkoMaldeghem01} characterizes twin buildings using the concept of $1$\ndash twinnings developed by Bernhard M\'ehlherr:

\begin{definition}[$1$\ndash twinning]
Let $\mathcal{C}^+$ (resp.\ $\mathcal{C}^-$) denote the chambers of $\mathfrak{B}^+$ (resp.\ $\mathfrak{B}^-$). A nonempty symmetric relation $\mathcal{O} \subset(\mathcal{C}_{+} \times \mathcal{C}_{-})\cup(\mathcal{C}_{-} \times \mathcal{C}_{+})$ is called a $1$\ndash twinning if the following axiom holds:
For a pair of chambers $(c_+, c_-) \in \mathcal{O}$ and two panels $P_+ \subset \mathcal{C_+}$ and $P_- \subset \mathcal{C_-}$ of the same type such that $c_{\epsilon} \in P_{\epsilon}, \epsilon \in \{+,-\}$ and $x_{\epsilon} \in P_{\epsilon}$ there exists a unique $y_{-\epsilon} \in P_{-\epsilon}$ such that $(x_{\epsilon}, y_{-\epsilon})\notin \mathcal{O}$.
\end{definition}

\begin{theorem}[Criterion for a twin building]
\label{criterionforatwinbuilding}
A $1$\ndash twinning induces a twin building iff for $\epsilon \in \{+,-\}$ there exists a chamber $c_{-\epsilon} \in \mathcal{C}_{-\epsilon}$ such that for any chamber $x_{\epsilon}$ with $(c_{-\epsilon}, x_{\epsilon})\subset \mathcal{O}$ there is an apartment $\Sigma_{\epsilon}$ of $\mathfrak{B}_{\epsilon}$ satisfying $\{x_{\epsilon}\}= \{y_{\epsilon}\in \mathcal{C}_{\epsilon}|y_{\epsilon} \in \Sigma_{\epsilon} \textrm{ and } (c_{-\epsilon}, y_{\epsilon} ) \in \mathcal{O}\}$.
\end{theorem}

The second characterization due to Peter Abramenko and Mark Ronan characterizes twin buildings via twin apartments~\cite{AbramenkoRonan98}. In complete analogy to the finite dimensional situation, twin apartments are characterized as the coconvex hull of opposite chambers. We quote their definition, adapting the notation to our conventions:

\begin{definition}[Twin apartment]
For every pair of opposite chambers $c_+$ and $c_-$, the set defined by $\{d_{\epsilon} \in \mathcal{C}^{\epsilon}| \delta^*(c_{-\epsilon}, d_{\epsilon})=\delta_{\epsilon}(c_{\epsilon}, d_{\epsilon})\}$ is an apartment in $\mathfrak{B}^{\epsilon}$ $(\epsilon \in \{+,-\})$. This apartment is called a twin apartment in $\mathfrak{B}$.
\end{definition}

\bibliographystyle{alpha}
\bibliography{Doktorarbeit}

\end{document}